\documentclass[UTF8,a4paper,12pt]{article}

\usepackage{amsthm}
\usepackage{natbib}
\usepackage{amsmath}
\usepackage{amssymb}
\usepackage{mathrsfs}
\usepackage{comment}
\usepackage{enumerate}
\usepackage{geometry}
\usepackage{mathabx}
\usepackage{hyperref}
\usepackage{dsfont}
\usepackage{microtype}

\usepackage{lipsum}
\usepackage{amsfonts}
\usepackage{graphicx}
\usepackage{epstopdf}
\usepackage{algorithmic}

\usepackage{enumitem}
\setlist[enumerate]{leftmargin=.5in}
\setlist[itemize]{leftmargin=.5in}

\theoremstyle{plain}
\newtheorem{theorem}{Theorem}
\newtheorem{corollary}{Corollary}
\newtheorem{lemma}{Lemma}
\newtheorem{proposition}{Proposition}

\theoremstyle{definition}
\newtheorem{definition}{Definition}

\newtheorem{assumption}{Assumption}

\theoremstyle{remark}

\newtheorem{problem}{Problem}

\newcommand{\ud}{\mathrm{d}}

\newcommand{\E}{\mathbb{E}}

\newcommand{\F}{\mathcal{F}}

\newcommand{\R}{\mathbb{R}}

\newcommand{\p}{\mathbb{P}}

\renewcommand{\(}{\left(}
\renewcommand{\)}{\right)}
\renewcommand{\[}{\left[}
\renewcommand{\]}{\right]}

\renewcommand{\geq}{\geqslant}
\renewcommand{\leq}{\leqslant}
\renewcommand{\epsilon}{\varepsilon}
\newcommand{\esssup}{\mathrm{ess\mbox{-}sup}}

\title{Mean-Field Game of Relative Performance Portfolio for Two Populations with Poisson Common Noise
}

\author{
	Yuchen Li\thanks{Department of Mathematical Sciences, Tsinghua University, Beijing 100084, China. Email: li-yc21@mails.tsinghua.edu.cn}
	\and
	Zongxia Liang\thanks{Department of Mathematical Sciences, Tsinghua University, Beijing 100084, China. Email: liangzongxia@tsinghua.edu.cn}
	\and
	Xiang Yu\thanks{Department of Applied Mathematics, The Hong Kong Polytechnic University, Kowloon, Hong Kong. Email: xiang.yu@polyu.edu.hk}
}

\begin{document}

\maketitle

\begin{abstract}
This paper studies the mean field game (MFG) and N-player game on relative performance portfolio management with two heterogeneous populations. In addition to the Brownian idiosyncratic and common noise, the first population invests in assets driven by idiosyncratic Poisson jump risk, while the second population invests in assets subject to Poisson common noise. We establish the characterization of the mean-field equilibrium (MFE) in MFG with two populations as well as the Nash equilibrium in the $N_1+N_2$-player game. Furthermore, we prove the convergence of the Nash equilibrium in the $N_1+N_2$-player game to the MFE as the number of players in two populations tends to infinity. We also discuss some impacts on MFE by the Poisson idiosyncratic risk and Poisson common noise in the context of relative performance, compensated by some numerical examples and financial implications.
\end{abstract}

\newcommand\keywords[1]{\textbf{Keywords}: #1}



\section{Introduction}
Initially introduced by {Lasry and Lions \cite{LL2007Mean}} and {Huang et al. \cite{HCM2007Large}}, mean-field game (MFG) has become an important and rapidly growing paradigm for analyzing strategic decisions for large population systems with interacting agents. {Lasry and Lions \cite{LL2007Mean}} established the analytical method for MFG by analyzing the coupled Hamilton-Jacobi-Bellman (HJB) and Fokker-Planck equations. {Huang et al. \cite{HCM2007Large}} approached MFG from an engineering control perspective, studying linear-quadratic-Gaussian problems with cost-coupled large populations and proposing a decentralized $\epsilon$-Nash equilibrium. These pioneer studies laid the foundation for recent advancements in MFG in both theories and applications.
	
	For fund management in a competition setting, relative performance criteria have attracted an upsurge of interest recently, which focus on how each agent optimizes the portfolio concerning other peers' portfolio performance. Typically, each agent sets the average performance of the entire population as a benchmark in her objective functional. {Espinosa \cite{Espinosa2015relative}} studied the stochastic differential game on relative performance portfolio management with exponential utility and trading constraints. {Lacker and Zariphopoulou \cite{LD2019mean}} revisited the relative performance game with CARA and CRRA utilities under asset specialization, in which they derived explicit constant equilibrium strategies in both finite-player and mean-field games. This study spurred a considerable amount of subsequent research on the topic of relative performance MFG in various context. {Lacker and Soret \cite{lacker_many-player_2020}} studied a similar relative performance MFG by concerning both terminal wealth and intermediate consumption. {Kraft et al. \cite{KRAFT2020103857}} discussed the portfolio optimization problems of two agents with relative wealth concerns and stochastic volatility. Furthermore, {dos Reis and Platonov \cite{DP2021forward}} explored the relative performance problem within the forward utility framework. Fu and Zhou \cite{fu_mean_2023} utilized the BSDE approach to investigate the relative performance problem with random coefficients. {Wang and Hu \cite{Wang04032023}} considered a family of portfolio management problems under relative performance competition and the memory effect. {Tangpi and Zhou \cite{tangpi_optimal_2024}} examined a class of graphon MFG under relative performance portflio. {Souganidis and Zariphopoulou \cite{souganidis_mean_2024}} investigated the master equation of MFG under relative performance under general utility. {Liang and Zhang \cite{LZinconsistent} }investigated a class of time-inconsistent MFG in the relative performance framework. {Bo et al. \cite{BWY2024Habit}} and {Liang and Zhang \cite{LZ2024habit}} applied MFG framework to study the relative consumption problem in a large population to incorporate the external habit formation from peer agents.
	
	
	
	 On the other hand, Poisson jump risk has caught increasing attention in mean-field game problems in different context. {Benazzoli \cite{BCD2020jump}} considered a general MFG problem with the idiosyncratic Poisson jump process using the compactification and relaxed control formulation. {Bo et al. \cite{BWY2024Mean}} extended the relative performance problem to include the idiosyncratic contagious jump risk modeled by multidimensional nonlinear Hawkes processes. {Hernández-Hernández and Ricalde-Guerrero \cite{HH2024poisson}}  studied general MFG with Poisson common noise using the maximum principle and BSDE approach. {Bo et al. \cite{BWWY}} recently established the existence of strong mean field equilibrium in MFG with Poissonian common noise using the relaxed control and pathwise formulation approach.
    The work that is most similar to ours is {Bank and Sedrakjan\cite{PG2025careknowjumpsignals}}, in which they considered Poisson common noise and a type of relative performance MFG using the PDE approach. However, they considered a special problem formulation where the investor is assumed to have future information of the Poisson risk and chooses the control that depends solely on the jump risk.

    In contrast to the existing literature, we consider two heterogeneous populations in the relative performance MFG and aim to analyze distinct behavior of investors when the Poisson jump risk appears in one population as the idiosyncratic noise and in another population as the common noise. Each population concerns the weighted average of their own population performance as well as the performance by the other population. This average is depicted by the weighted conditional mean with respect to the common noise filtration. Such multi-population games are relatively new in the context of portfolio management under the relative performance criteria. First, we analyze MFG problem between two heterogeneous populations using the PDE approach and some fixed point argument. Then, we examine the $N_1+N_2$-player game involving $N_1$ investors in population 1 and $N_2$ investors in population 2. We derive both the mean-field equilibrium in the MFG and the Nash equilibrium in the $N_1+N_2$-player game for representative agents in two populations, explicitly illustrating distinct equilibrium portfolio strategy when the Poisson jump only appers as the idiosyncratic noise (in population 1) and plays the role of common noise (in populations 2). Furthermore, we also rigorously prove the limit theory to show that the mean-field equilibrium is the limit of the Nash equilibrium strategy in the $N_1+N_2$-player game as $N_1,N_2$ tend to infinity. 
	
	In both the MFG and the $N_1+N_2$-player game, we adopt an approach similar to that in {Lacker and Zariphopoulou \cite{LD2019mean}}. First, we consider the auxiliary problems for representative agents in two populations given the filtration generated by their common noise, in which the HJB equations are derived and the verification theorems for best response controls are established. Second, we characterize the equilibrium as the fixed point under the consistency condition. In both problems, the existence and uniqueness of the equilibrium reduce to the study of a three-dimensional system of equations. It is interesting to observe that the Poisson jump risk distorts the linear structure of the equilibrium strategy derived in {Lacker and Zariphopoulou \cite{LD2019mean}}. Investors facing Poisson jump risk exhibit nonlinear interactions with relative performance. We discuss some quantitative impacts of incorporating Poisson jump risk and present some numerical examples to compensate the theoretical findings and discuss their financial implications. 
	
	The rest of the paper is organized as follows. In Section \ref{sec: formulation}, we introduce the MFG problem with two populations such that the Poisson jump risk appears as the idiosyncratic noise and common noise respectively, in which the mean field equilibrium strategies for different populations are characterized. In Section \ref{sec: N1N2player}, we define and solve the associated $N_1+N_2$-player game. In Section \ref{sec: convergence}, we rigorously prove the convergence from the Nash equilibrium from the $N_1+N_2$ player game to the mean field equilibrium as $N_1,N_2$ tend to infinity. In Section \ref{sec: numerical}, we analyze the influence of Poisson common noise and present several numerical examples of sensitivity analysis with financial discussions.

	\section{Problem Formulation}\label{sec: formulation}
	
    We consider a filtered probability space $(\Omega;\F;\{\F_t\}_{t\geq 0};\p)$ satisfying the usual conditions. $W^i=\{W^i_t,t\in[0,T]\},\ i=0,1,2$  are Brownian motions. ${N^1}$, ${N}$ are  Poisson random measures on $\([0,T]\times \R,\mathcal{B}([0,T]\times \R)\)$, respectively. Filtration $\{\F_t\}_{t\geq 0}$ is generated by $W^0$, $W^1$, $W^2$, $N^1$ and $N$. We assume that $W^i,\ i=0,1,2$, $N^1$ and $N$ are independent. Denote  $\ud t\nu^1(\ud p)$ and $\ud t\nu(\ud p)$ as the compensators (or intensity measures), and $\tilde{N}^1$ and $\tilde{N}$ as the corresponding compensated martingale measures, i.e., 
	$$\tilde{N^1}\(\ud t,\ud p\)=N^1\(\ud t,\ud p\)-\ud t\nu^1(\ud p),\ \tilde{N}\(\ud t,\ud p\)=N\(\ud t,\ud p\)-\ud t\nu(\ud p).$$
	It is assumed that $\nu^1(\R)<\infty, \nu(\R)< \infty$. 
	 Furthermore, let us assume that the filtration supports two random vectors, representing the information of two populations that
	$$\zeta^1=(\xi^1,\delta^1,\lambda^1,\mu^1,\sigma^1,\sigma^{0,1},\gamma^1),\ \zeta^2=(\xi^2,\delta^2,\lambda^2,\mu^2,\sigma^2,\sigma^{0,2},\gamma^2).$$
	In both random vectors, $\xi^i$, $\delta^i$ and $\lambda^i=(\lambda_{i,1},\lambda_{i,2} )$ are random variables represent the initial value, risk aversion coefficient and the relative competition coefficient of the population. 
	$\sigma^1$, $\sigma^2$, $\sigma^{0,1}$, $\sigma^{0,2}$, $\gamma^1$ and $\gamma^2$ are random volatility coefficients.
	
	The following assumption is imposed.
	\begin{assumption}\label{ass bound}		
		$\zeta^1$ and $\zeta^2$ are bounded random vectors.
	\end{assumption}
	In the mean field model, the wealth processes of the two populations are denoted by $X^1$ and $X^2$ respectively, with dynamics
	\begin{equation}\label{eq: MFG dynamic of X}
		\left\{
		\begin{aligned}
			\ud X^1_t&=\pi^1_t\mu^1\ud t+\pi^1_t\sigma^1\ud W^1_t+\pi^1_t\sigma^{0,1}\ud W_t^0+\pi^1_t\gamma^1\ud \tilde{N}_t^1,\quad X_0^1= \xi^1,\\
			\ud X^2_t&=\pi^2_t\mu^2\ud t+\pi^2_t\sigma^2\ud W^2_t+\pi^2_t\sigma^{0,2}\ud W_t^0+\pi^2_t\gamma^2\ud \tilde{N}_t,\quad  X_0^2= \xi^2.			
		\end{aligned}\right.
	\end{equation}	
	In population 1, $W^1$ and $\tilde{N}^1$ stand for idiosyncratic Brownian and Poisson noise; while in population 2, $W^2$ is the idiosyncratic Brownian noise and $\tilde{N}$ represents the Poisson common noise attached to population 2. Moreover, $W^0$ is the Brownian common noise that is affecting two populations. The common noise filtration $\F_0$ is generated by the common noise $W^0$ and $\tilde{N}$. Filtration $\F_1$ is generated by $W^1$, $W^0$ and $\tilde{N}^1$, and filtration $\F_2$ is generated by $W^2$, $W^0$ and $\tilde{N}$.
	\begin{definition}
		For $i=1,2$, we call a strategy $\pi^i$ admissible or $\pi\in \mathcal{A}^i$ if and only if 
		\begin{enumerate}
			\item $\pi^i$ is predictable with respect to $\F^i$.
			\item  The SDE \eqref{eq: MFG dynamic of X} under $\pi^i$ has a unique solution satisfying \ 
			$\E[\sup\limits_{0\leq s\leq T}|X^{i}_s|^2]<\infty.$ 
		\end{enumerate}
	\end{definition}
In our MFG setting, both populations will consider an optimal investment problem with their own relative points. Suppose that the investor is able to see the conditional expectation of the final wealth of both populations $\E\[X^{i}_T|\F^0_T\]$, and $\lambda$ represents the relative competition coefficient of the two populations. Then, for population $i$, the relative point $H^i$ is $$H^i=\lambda_{i,1}\E\[X^{1}_T|\F^0_T\]+\lambda_{i,2}\E\[X^{2}_T|\F^0_T\].$$ 
 We consider the following relative performance MFG problem under exponential utility.
	\begin{problem}\label{prob: invest}
		\begin{equation}
			\max_{\pi^i} \E\Big[-\frac{1}{\delta^i}e^{-\delta^i \(X^{i}_T-H^i\)}\Big].
		\end{equation}
	\end{problem}
We next give the definition of the equilibrium strategy for our MFG problem with common noise $W^0$ and $\tilde{N}$: 
	\begin{definition}
		We say that $(\pi^{1,*},\pi^{2,*})$ is a Mean Field Equilibrium (MFE) if and only if $(\pi^{1,*},\pi^{2,*})$ is the optimal strategy of Problem \ref{prob: invest} with respect to relative point $$H^i=\lambda_{i,1}\E\[X^{1,*}_T|\F^0_T\]+ \lambda_{i,2}\E\[X^{2,*}_T|\F^0_T\].$$
	\end{definition}
	First, we are going to find an equilibrium strategy as an function of the parameters $\zeta^1$ and $\zeta^2$.
	
	\subsection{Finding equilibrium}
	%
	To derive the characterization of the MFE in two pupolations, we follow the standard procedure:
	\begin{enumerate}
		\item Fix $H^1,H^2\in \F_T$  and solve the best response control $\pi^{i,*}$, $i=1,2$, for the representative agent in each population in Problem \eqref{prob: invest}: $$
		\max_{\pi^1} \E\Big[-\frac{1}{\delta^1}e^{-\delta^1 \(X^{1}_T-H^1\)}\Big],\ 
		\max_{\pi^2} \E\Big[-\frac{1}{\delta^2}e^{-\delta^2 \(X^{2}_T-H^2\)}\Big]. $$
		\item Verify the consistency condition by solving the fixed point problem: 
		$$\begin{aligned}
			H^1=\lambda_{1,1}\E\[X^{1,*}_T|\F^0_T\]+\lambda_{1,2}\E\[X^{2,*}_T|\F^0_T\],\\
			H^2=\lambda_{2,1}\E\[X^{1,*}_T|\F^0_T\]+\lambda_{2,2}\E\[X^{2,*}_T|\F^0_T\],
		\end{aligned}$$where  $X_T^{1*}$  and $X_T^{2*}$  is the optimal wealth of step 1.
	\end{enumerate}
    \vskip 4pt
	The key step is to solve Problem \ref{prob: invest} first, for every proper relative point $H^i\in \F^0_T$. Similar to {Lacker and Zariphopoulou \cite{LD2019mean}}, for any $\F_T^0$ measurable random variables $\bar{X}_1$ and $\bar{X}_2$,  let $H^i=\lambda^i_1\bar{X}_1+\lambda^i_2\bar{X}_2$, then we have
	\begin{equation}
		\max_{\pi\in \mathcal{A}} \E\Big[-\frac{1}{\delta^i}e^{-\delta^i \(X^{i}_T-\(\lambda^i_1\bar{X}_1+\lambda^i_2\bar{X}_2\)\)}\Big]=\E\[V^i(\xi^i,0)\],
	\end{equation}
	where $V^i(\cdot)$ is the value function defined by the Problem \ref{prob: invest}, with deterministic type vector $\zeta^i$. 
    \vskip 5pt
    We restate this deterministic version of Problem \ref{prob: invest} and value function $V^i$ as follows. For simplicity and unity of the notations with respect to two different populations, we can assume that population 1 has $\gamma^2\equiv0$, population 2 has $\gamma^1\equiv0$, and ignore the subscript $i$ in the notations.
	\begin{problem}\label{prob: deter invest prob}
		Consider a specific augmented deterministic sample vector $$\zeta_0=(x,\delta,\lambda,\mu,\sigma,\sigma^0,\gamma^1,\gamma)$$ drawn from the random type vector $(\zeta,\delta,\lambda,\mu,\sigma,\sigma^0,\gamma^1,\gamma)$. Consider the dynamics of wealth process 
		\begin{equation}\label{eq: dynamic of determin X}
			\begin{aligned}
				\ud X_t=\pi_t\mu_t\ud t+\pi_t\sigma\ud W_t+\pi_t\sigma^0\ud W_t^0+\pi_t\gamma_t\ud \tilde{N}_t+\pi_t\gamma^1\ud \tilde{N}^1_t,\quad
				X_0=x,
			\end{aligned}
		\end{equation} and admissible strategy $\pi$ and admissible set $\mathcal{A}_0$ defined by
		\begin{enumerate}
			\item $\pi$ is predictable with respect to $\F$.
			\item The SDE \eqref{eq: dynamic of determin X} under $\pi$ has a unique solution satisfying
			$$\E\Big[\sup_{0\leq s\leq T}|X^{\pi}_s|^2\Big]<\infty.$$
		\end{enumerate}
		Then, the optimal relative performance problem with value function $V$ and relative potint $H=\lambda_1\bar{X}_1+\lambda_2\bar{X}_2$ is defined by 
		\begin{eqnarray}\label{probm2}
			V(x,t)=\max_{\pi\in \mathcal{A}_0} \E_{x,t}\Big[-\frac{1}{\delta}e^{-\delta \(X^\pi_T-\(\lambda_1\bar{X}_1+\lambda_2\bar{X}_2\)\)}\Big].
		\end{eqnarray}
	\end{problem}
    
	After solving the optimal strategy under deterministic problem, we can represent the solution $\pi^*(\zeta)$ of the Problem \ref{prob: invest} using $\pi^*(\zeta_0)$ from above deterministic problem by plugging in the random vector $\zeta$ into $\pi^*(\cdot)$. The requirement for the admissible and optimal condition for  $\pi^*(\zeta)$ is automatically satisfied. 
	\subsection{Solution of the auxiliary Problem \ref{prob: deter invest prob}}
	In this subsection, we will find the solution of the  problem \ref{prob: deter invest prob}. It is assumed that $H$ admits the dynamics with deterministic coefficients $\eta=\{\eta_t\}_{0\leq t\leq T}$, $u=\{u_t\}_{0\leq t\leq T}$ and $v=\{v_t\}_{0\leq t\leq T}$:
	$$\ud H_t=\eta_t\ud t+u_t\ud W_t^0+v_t\ud \tilde{N}_t.$$ 
	The meaning of $u$ and $v$ are the risk loads of the relative point process based on Brownian common noise and Poisson common noise.
	The difference between the representative agent's wealth and relative performance is $Z^\pi=X-H$ satisfying
	$$
	\begin{aligned}
		\ud Z_t^\pi=(\pi_t\mu_t-\eta_t)\ud t+\pi_t\sigma\ud W_t+(\pi_t\sigma^0-u_t)\ud W_t^0+(\pi_t\gamma_t-v_t)\ud \tilde{N}_t+\pi_t\gamma^1\ud \tilde{N}^1_t.
	\end{aligned}
	$$
	Equivalently, the deterministic version of \ref{prob: deter invest prob}  is given by the next auxiliary problem:  
	\begin{problem}\label{prob MFE determin Z}
		\begin{eqnarray}
			V(x,t)=\max_{\pi\in\mathcal{A}_0}\E\Big[-\frac{1}{\delta}e^{-\delta Z_T^\pi}|Z_t=x\Big].
		\end{eqnarray}
	\end{problem}
	According to {{\O}ksendal and Sulem \cite{Oksendal2007ch3}}, the generator $\mathscr{A}$ of $Z$ is
	\begin{equation}\label{MFG generator}
		\begin{aligned}
			\mathscr{A}^{\pi} V(x,t)=&V_t+V_x\[\pi_t\mu_t-\eta_t-\(\pi_t\gamma_t-v_t\)\nu-\pi_t\gamma^1\nu^1\]
			+\frac{1}{2}V_{xx}\[\pi_t^2\sigma^2+\(\pi_t\sigma^0-u_t\)^2\]\\
			&+\[V(x+\pi_t\gamma_t-v_t,t)-V(x,t)\]\nu+\[V(x+\pi_t\gamma^1,t)-V(x,t)\]\nu^1.
		\end{aligned}
	\end{equation}
	
	To solve Problem \ref{prob MFE determin Z}, we will prove the following verification theorem. For simplicity, we slightly abuse the notation $\pi$ to represent a process or a real number in $\R$. Denote $J^{\pi}(x,t)=\E\big[-\frac{1}{\delta}e^{-\delta Z_T^\pi}|Z_t=x\big]$.
	
	\begin{theorem}\label{thm: verification}
		Suppose that there exist a function $  V\in C^{2,1}\(\R\times[0,T)\)\cap C\(\R\times[0,T]\)$ and a control $ \pi^*\in \mathcal{A}_{0}$  such that 
		\begin{enumerate}[label=(\roman*)]
			\item $\mathscr{A}^{\pi^*}V(x,t)= 0,\ \forall \(x,t\) \in \R \times [0,T);$
			\item $\mathscr{A}^{\pi}V(x,t)\leq 0,\ \forall \(x,t\) \in \R \times \R^+ \times \R \times [0,T),\ \forall \pi \in \R;$
			\item $V(x,T)=-\frac{1}{\delta}e^{-\delta x}\ a.s.,\  \forall x \in \R;$
			\item $\E_{x,t}\Big[\sup\limits_{t\leq s\leq T}|V(Z_s^{\pi},s))|\Big]<+\infty,\ \text{for} \ \forall\(x,t\) \in \R \times [0,T), \forall \pi \in \mathcal{A}_{0}$.
		\end{enumerate}
		Then 
		\begin{eqnarray*}
			&&J^{\pi}(x,t)\leq V(x,t),\ \text{for} \ \forall \pi\in\mathcal{A}_{0},\\
			&&V(x,t)=J^{\pi^*}(x,t)=\max_{\pi} J^{\pi}(x,t),
		\end{eqnarray*}
		and $\pi^*$ is the optimal control and $V(x,t)$ is the value function of Problem \ref{prob MFE determin Z}. 
	\end{theorem}
	
	\begin{proof}
		For $\forall \pi\in \mathcal{A}_{0}$, define a stopping time $\tau_N=T\wedge \inf\{t\geq 0: |Z_t|\geq N\}$. Using  Ito's formula (see, e.g., {Protter \cite{Protter2005}}), we have 
		\begin{equation}\label{dynkin formula}
			\E_{x,t}\[V(Z_{\tau_N},\tau_N)\]=V(x,t)+\E_{x,t}\Big[\int_{t}^{\tau_N} \mathscr{A}^{\pi}V(Z_{s-},s)\ud s\Big].
		\end{equation}
		First,  taking $\pi^*$ into Equation \eqref{dynkin formula} and using Condition $(\romannumeral1)$, we have $$\E_{x,t}\[V(Z_{\tau_N}^{\pi^*},\tau_N)\]= V(x,t).$$ 
		Using Conditions $(\romannumeral3)$, $(\romannumeral4)$ and applying Dominated Convergence Theorem as $N\to \infty$, we obtain $$J^{\pi^*}(x,t)= \E_{x,t}\Big[-\frac{1}{\delta}e^{-\delta Z_T^{\pi^*}}\Big]= V(x,t).$$ 
		
		Second, Equation \eqref{dynkin formula} and Condition $(\romannumeral2)$ yield $$\E_{x,t}\[V(Z_{\tau_N},\tau_N)\]\leq V(x,t).$$ Letting $N\to \infty$, we obtain $$J^{\pi}(x,t)\leq V(x,t),$$
		which completes the proof.
	\end{proof}
	We make the ansatz of the solution   
    to Problem \ref{prob MFE determin Z}. Based on the boundary condition $(\romannumeral3)$, we conjecture that $V(x,t)=-\frac{1}{\delta}e^{-\delta x}f(t)$. By virtue of Condition $(\romannumeral1)$, $f$ should satisfy 
	\begin{equation}\label{eq: ODE f}
		\begin{aligned}
			0=&f'(t)+f(t)\min_{\pi_t}\Big\{-\delta\Big[\pi_t\mu_t-\eta_t-(\pi_t\gamma_t-v_t)\nu-\pi_t\gamma^1\Big]\\&+\frac{1}{2}\delta^2\Big[\pi_t^2\sigma^2+(\pi_t\sigma^0-u_t)^2\Big]+\Big[e^{-\delta(\pi_t\gamma_t-v_t)}-1\Big]\nu+\Big[e^{-\delta(\pi_t\gamma^1)}-1\Big]\nu^1\Big\}
		\end{aligned}
	\end{equation}
	with $f(T)=1$.
	
	Based on Equation \eqref{eq: ODE f}, it follows  that
	$$\pi^*_t=\arg\min_{\pi} F(\pi,\alpha,t),$$
	where $F$ is defined by
	$$F(\pi,t)=\frac{1}{2}\delta^2\[\pi^2\sigma^2+\(\pi\sigma^0-u_t\)^2\]+\[e^{-\delta\(\pi\gamma_t-v_t\)}+\delta\pi\gamma_t\]\nu+\[e^{-\delta\(\pi\gamma^1\)}+\delta\pi\gamma^1\]\nu^1-\delta\pi\mu_t.$$
	Let $g(\pi,t)=\frac{\partial F}{\partial \pi}(\pi,t)$. Then $$g(\pi,t)=-\delta \gamma_t\nu \[e^{-\delta\(\pi\gamma_t-v_t\)}-1\]-\delta \gamma^1\nu^1 \[e^{-\delta\(\pi\gamma^1\)}-1\]+\delta^2\(\sigma^2+{\sigma^0}^2\)\pi-\delta^2{\sigma^0}u_t-\delta\mu_t.$$
	We see that $\lim\limits_{\pi\to -\infty}g(\pi,t)=-\infty$, $\lim\limits_{\pi\to +\infty}g(\pi,t)=+\infty$ and $g$ is a strictly monotone with respect to $\pi$. Thus, for every $(\omega,t)\in \Omega\times [0,T]$, there exists a unique $\pi$ such that $g(\pi,t)=0$. We denote the root of $g(\pi,t)$ as $\pi^*$ and aim to verify that $\pi^*$ is the optimal control for Problem \ref{prob MFE determin Z}.
	\begin{theorem}\label{thm: MFG optimal pi}
		Let $\pi^*$ be the unique solution of the equation $g(\pi,t)=0$, and $V(x,t)=-\frac{1}{\delta}e^{-\delta x}f(t)$, where $f$ satisfies the ODE
		\begin{equation*}\left\{
			\begin{aligned}
				0=&f'(t)+f(t)\Big\{-\delta\[\pi^*_t\mu_t-\eta_t-\(\pi^*_t\gamma_t-v_t\)\nu-\(\pi^*_t\gamma^1\)\nu^1\]\\&+\frac{1}{2}\delta^2\[\pi_t^2\sigma^2+\(\pi^*_t\sigma^0-u_t\)^2\]+\[e^{-\delta\(\pi^*_t\gamma_t-v_t\)}-1\]\nu+\[e^{-\delta\(\pi^*_t\gamma^1\)}-1\]\nu^1\Big\},\\
				1=&f(T).
			\end{aligned}\right.
		\end{equation*}
		Then $\pi^*$ and $V(x,t)$ are the best response control and the associated value function of Problem \ref{prob MFE determin Z}.
	\end{theorem}
	\begin{proof}
		We only need to show that $\pi^*_t$ and $V(x,t)$ satisfy the verification theorem \ref{thm: verification}. Conditions $(\romannumeral1)$,  $(\romannumeral2)$ and $(\romannumeral3)$ are naturally satisfied by the definition of $f$.
		As we can see that $\pi^*_t$ is a continuous function of our market  coefficients,  $u_t$ and $v_t$. Thus, $\pi^*_t$ and $f(t)$ are bounded according to our Assumption \ref{ass bound}, which implies that $V(Z^{\pi^*}_t,t)$ satisfies a linear SDE with bounded coefficient, and Condition $(\romannumeral4)$ holds.
	\end{proof}
	Now, consider Problem \ref{prob MFE determin Z} for the population 1 and population 2 respectively, we can derive the distinct characterization of their best response controls. For the population 1, we have $\sigma=\sigma^1$, $W=W^1$ and $\gamma=0$. For the population 2, we have $\sigma=\sigma^2$, $W=W^2$ and $\gamma^1=0$.
	\begin{proposition}\label{prop: MFG optimal pi}
		For two populations 1 and 2, their best response control $\pi^{1,*}$ and $\pi^{2,*}$ are given by the unique solution of the following two equations,respectively:
		$$g^1(\pi,t)=-\delta \gamma^1\nu^1 \[e^{-\delta\(\pi\gamma^1\)}-1\]+\delta^2\({\sigma^1}^2+{\sigma^{0,1}}^2\)\pi-\delta^2{\sigma^{0,1}}u_t-\delta\mu^1,
        $$
		$$g^2(\pi,t)=-\delta \gamma\nu \[e^{-\delta\(\pi\gamma-v_t\)}-1\]+\delta^2\({\sigma^2}^2+{\sigma^{0,2}}^2\)\pi-\delta^2{\sigma^{0,2}}u_t-\delta\mu^2.$$
	\end{proposition}
	The proposition above gives the optimal strategy under fixed $u,v$. As we can see that $g^1$ is independent of $v$, which indicates that the best response strategy of the  population 1 will not be influenced by the Poisson common noise, which is an expected result. 
	
	In the simple case as in {Lacker and Zariphopoulou \cite{LD2019mean}}, for the population 1, if there is no idiosyncratic Poisson noise, we have the next result.
	\begin{proposition}\label{prop: MFG optimal pi 1}
		For the population 1, if $\gamma^1=0$, the optimal best response control is $$\pi^{1,*}_t=\frac{{\sigma^{0,1}}}{{\sigma^1}^2+{\sigma^{0,1}}^2}u_t+\frac{\mu_t}{\delta(\sigma^2+{\sigma^0}^2)}.$$
	\end{proposition} 
	Note that $f,g^1,g^2,\pi^*$ are in fact functions of market coefficients,  $u_t$ and $v_t$. Therefore, from now on, we denote
	$$g(\pi,u,v,\zeta)=-\delta \gamma\nu \[e^{-\delta\(\pi\gamma- v\)}-1\]+\delta^2\(\sigma^2+{\sigma^0}^2\)\pi-\delta^2{\sigma^0}u-\delta\mu.$$
	It follows that $$g^1(\pi,t)=g(\pi,u,0,\zeta^1),\ g^2(\pi,t)=g(\pi,u,v,\zeta^2).$$  The optimal strategy can be regarded as a function of $u$ and $v$: $\pi^{1,*}(u,v)$ is the unique root of equation $g(\pi,u,0,\zeta^1)=0$, $\pi^{2,*}(u,v)$ is the unique root of equation $g(\pi,u,v,\zeta^2)=0$. 

	\subsection{Solution of the MFG}\label{MFG}
	In this subsection, we will characterize the MFE for each population. To simplify notations, we omit the time subscript $t$ whenever no confusion arises.
	
	First, assume that the deterministic equilibrium strategy exists and denote it by $(\alpha^1,\alpha^2)=\{(\alpha^1_t,\alpha^2_t)\}_{0\leq t\leq T}$.  Then, the dynamics of the mean wealth processes  are
	$$\begin{aligned}
		&\ud \widebar{X^{1,\alpha^1}_t}=\widebar{\alpha^1_t\mu^1}\ud t+\widebar{\alpha^1_t\sigma^{0,1}}\ud W_t^0,\\
		&\ud \widebar{X^{2,\alpha^2}_t}=\widebar{\alpha^2_t\mu^2}\ud t+\widebar{\alpha^2_t\sigma^{0,2}}\ud W_t^0+\widebar{\alpha^2_t\gamma}\ud \tilde{N}_t,
	\end{aligned}$$
	 where $\bar{\cdot}$ stands for the conditional expectation given the filtration $\F^0$ (e.g., $\widebar{\alpha^1_t\mu^1}=\E\[\alpha^1_t\mu^1|\F^0\]$). Because $\alpha$, $\mu$, $\sigma^0$ and $\gamma$ are independent of $\F_T^0$,  the $\widebar{\alpha_t\mu}$, $\widebar{\alpha_t\sigma^0}$ and $\widebar{\alpha_t\gamma}$ degenerate to deterministic functions $\E\[{\alpha_t\mu}\]$, $\E\[{\alpha_t\sigma^0}\]$ and $\E\[{\alpha_t\gamma}\]$.
	Now, with the notations: $x_{1,t}=\E\[{\alpha^1_t\sigma^{0,1}}\]$, $x_{2,t}=\E\[{\alpha^2_t\sigma^{0,1}}\]$ and  $y_t=\E\[{\alpha^2_t\gamma}\]$, the relative points $H^1$ and $H^2$ satisfy
	$$\begin{aligned}
		&\ud H^1_t=(\lambda_{1,1}\widebar{\alpha^1_t\mu^1}+\lambda_{1,2}\widebar{\alpha^2_t\mu^2})\ud t+(\lambda_{1,1}x^1_t+\lambda_{1,2}x^2_t)\ud W_t^0+\lambda_{1,2}y_t\ud \tilde{N}_t,\\ &\ud H^2_t=(\lambda_{2,1}\widebar{\alpha^1_t\mu^1}+\lambda_{2,2}\widebar{\alpha^2_t\mu^2})\ud t+(\lambda_{2,1}x^1_t+\lambda_{2,2}x^2_t)\ud W_t^0+\lambda_{2,2}y_t\ud \tilde{N}_t.
	\end{aligned}$$
	Thus, for the auxiliary control Problem \ref{prob: deter invest prob} of population 1, $\eta^1_t=\lambda_{1,1}\widebar{\alpha^1_t\mu^1}+\lambda_{1,2}\widebar{\alpha^2_t\mu^2}$, 
	$u^1_t=\lambda_{1,1}x^1_t+\lambda_{1,2}x^2_t$, $v^1_t=\lambda_{1,2}y_t$. For  population 2, $\eta^2_t=\lambda_{2,1}\widebar{\alpha^1_t\mu^1}+\lambda_{2,2}\widebar{\alpha^2_t\mu^2}$, 
	$u^2_t=\lambda_{2,1}x^1_t+\lambda_{2,2}x^2_t$, $v^2_t=\lambda_{2,2}y_t$.
	According to Theorem \ref{thm: MFG optimal pi} and Proposition \ref{prop: MFG optimal pi}, we get 
	$$\pi^{1,*}_t=\pi^{1,*}(u^1_t,v^1_t),\ \ 
	\pi^{2,*}_t=\pi^{2,*}(u^2_t,v^2_t),$$
	where $\pi^{1,*}(u,v)$ and $\pi^{2,*}(u,v)$ are the solutions of 
	$g(\pi,u,0,\zeta^1)=0$ and $g(\pi,u,v,\zeta^2)=0$, respectively.
	As $\pi^{1,*}(u,v)$ is independent of $v$, we can always assume that $v^1_t\equiv0$. 
	
	Consequently, we denote  $\pi^{*}(u,v)$ as the solution of 
	$g(\pi,u,v,\zeta)=0 $ with a univariate virtual type vector $\zeta$. 
	
	Given $\alpha$, we  find $\pi^{1,*}$ and $\pi^{2,*}$ with respect to the relative point driven by $\alpha$.
    
    If the consistency condition holds that $\alpha=(\pi^{1,*},\pi^{2,*})$, we find the desired equilibrium. On the other hand, we can easily find that the optimal strategies $\pi^{1,*}$ and $\pi^{2,*}$ only depend on unknown deterministic processes $x^1,x^2$ and $y$. Thus, to find a constant equilibrium, we could just find the fixed point of the following system of three equations: \begin{equation}\label{eq: two player eq}
		\left\{\begin{aligned}
			&\E\[\pi^{1,*}(\lambda_{11} x_1+\lambda_{12}x_2,0)\sigma^{0,1}\]=x_1,\\
			&\E\[\pi^{2,*}(\lambda_{21} x_1+\lambda_{22}x_2,\lambda_{22}y)\sigma^{0,2}\]=x_2,\\
			&\E\[\pi^{2,*}(\lambda_{21} x_1+\lambda_{22}x_2,\lambda_{22}y)\gamma\]=y.\\
		\end{aligned}\right.
	\end{equation}
	As long as Equation \eqref{eq: two player eq} has a fixed point $(x^*_1,x^*_2,y^*)$, we have an equilibrium $$(\pi^{1,*}(\lambda_{11} x^*_1+\lambda_{12}x^*_2,0),\pi^{2,*}(\lambda_{11} x^*_1+\lambda_{12}x^*_2,\lambda_{12}y^*)).$$
	By imposing some proper conditions on $\lambda$, we can establish the existence and uniqueness of solution to equation system \eqref{eq: two player eq}. The Implicit function theorem yields \begin{equation}\label{eq: MFG partial eq}
		\begin{aligned}
			\frac{\partial\pi^*(u,v)}{\partial u}=-\frac{\frac{\partial g}{\partial u}}{\frac{\partial g}{\partial \pi}}=\frac{\sigma^{0,2}}{({\sigma^{0,2}}^2+\sigma^2)+\gamma^2\nu e^{-\delta\(\pi\gamma- v\)}},\\
			\frac{\partial\pi^*(u,v)}{\partial v}=-\frac{\frac{\partial g}{\partial v}}{\frac{\partial g}{\partial \pi}}=\frac{\gamma\nu e^{-\delta\(\pi\gamma-\theta v\)}}{({\sigma^{0,2}}^2+\sigma^2)+\gamma^2\nu e^{-\delta\(\pi\gamma- v\)}}.
		\end{aligned}
	\end{equation}
	We prove that for some small $\lambda$,  equation system \eqref{eq: two player eq} has a unique solution.
	\begin{theorem}\label{thm: MFG existent of equi}
		There exists an $\epsilon>0$ depending on parameters $\sigma^{0,2}$ and $\gamma$  such that when $\|\lambda\|_{\infty}<\epsilon$, equation system \eqref{eq: two player eq} admits a unique solution $(x_1,x_2,y)\in \R^3$.
	\end{theorem}
	\begin{proof}
		We use the contracting mapping to draw the conclusion. Define $$R(x_1,x_2,y)=\(\begin{matrix}
			&\E\[\pi^{1,*}(\lambda_{11} x_1+\lambda_{12}x_2,0)\sigma^{0,1}\]\\
			&\E\[\pi^{2,*}(\lambda_{21} x_1+\lambda_{22}x_2,\lambda_{22}y)\sigma^{0,2}\]\\
			&\E\[\pi^{2,*}(\lambda_{21} x_1+\lambda_{22}x_2,\lambda_{22}y)\gamma\]\\
		\end{matrix}\).$$ We only need to prove that there exists some $L<1$ for  the inequality:
        $$\|R(x,y)-R(x',y')\|\leq L\|(x,y)-(x',y')\|.$$
		First, based on the Lagrange mean value theorem, we have 
		$$\begin{aligned}
			&\E\[\pi^{2,*}(\lambda_{21} x_1+\lambda_{22}x_2,\lambda_{22}y)\sigma^{0,2}\]-\E\[\pi^{2,*}(\lambda_{21} x'_1+\lambda_{22}x'_2,\lambda_{22}y')\sigma^{0,2}\]\\
			=&\E\Bigg[\sigma^{0,2}\[\lambda_{21}(x_1-x'_1)+\lambda_{22}(x_2-x'_2)\]\frac{\partial\pi^{2,*}}{\partial u}+\sigma^{0,2}\[\lambda_{22}(y-y')\]\frac{\partial\pi^{2,*}}{\partial v}\Bigg]\\
			\leq &\|(x_1,x_2,y)-(x'_1,x'_2,y')\|\E\[\sqrt{({\sigma^{0,2}})^2\Pi}\],
		\end{aligned}$$
        where $$\Pi:=\(\lambda_{21}\frac{\partial\pi^{2,*}}{\partial u}\)^2+\(\lambda_{22}\frac{\partial\pi^{2,*}}{\partial u}\)^2+\(\lambda_{22}\frac{\partial\pi^{2,*}}{\partial v}\)^2,$$
        evaluated at $$(u,v)=\lambda_{21}(x'_1+\theta(x_1-x'_1))+\lambda_{22}(x'_2+\theta(x_2-x'_2)),\lambda_{22}(y'+\theta(y-y')).$$
		Using Jensen's inequality, we have 
		\begin{eqnarray*}
			&&\(\E \[\pi^{2,*}(\lambda_{21} x_1 +\lambda_{22}x_2,\lambda_{22}y)\sigma^{0,2}\]-\E\[\pi^{2,*}(\lambda_{21} x'_1+\lambda_{22}x'_2,\lambda_{22}y')\sigma^{0,2}\]\)^2\\
			&&\leq\|(x_1,x_2,y)-(x'_1,x'_2,y')\|^2\E\[\sqrt{({\sigma^{0,2}})^2\Pi}\]^2\\
			&&\leq\|(x_1,x_2,y)-(x'_1,x'_2,y')\|^2\E\[({\sigma^{0,2}})^2\Pi\].
		\end{eqnarray*}
		Again for the exactly same $\theta$, we repeat the same process to get 
        \begin{eqnarray*}
			&&\(\E \[\pi^{2,*}(\lambda_{21} x_1+\lambda_{22}x_2,\lambda_{22}y)\gamma\]-\E\[\pi^{2,*}(\lambda_{21} x'_1+\lambda_{22}x'_2,\lambda_{22}y')\gamma\]\)^2\\
			&&\leq\|(x_1,x_2,y)-(x'_1,x'_2,y')\|^2\E\Big[{{\gamma}^2\big(\big(\lambda_{21}\frac{\partial\pi^{2,*}}{\partial u}\big)^2+\big(\lambda_{22}\frac{\partial\pi^{2,*}}{\partial u}\big)^2+\big(\lambda_{22}\frac{\partial\pi^{2,*}}{\partial v}\big)^2\big)}\Big],
		\end{eqnarray*}
		and 
        \begin{eqnarray*}
			&&\(\E\[\pi^{1,*}(\lambda_{21} x_1+\lambda_{22}x_2,0)\gamma\]-\E\[\pi^{1,*}(\lambda_{21} x'_1+\lambda_{22}x'_2,0)\gamma\]\)^2\\
			&&\leq\|(x_1,x_2,y)-(x'_1,x'_2,y')\|^2\E\Big[{({\sigma^{0,1})}^2\big(\lambda_{11}^2\big(\frac{\partial\pi^{1,*}}{\partial u}\big)^2+\lambda_{12}^2\big(\frac{\partial\pi^{1,*}}{\partial u}\big)^2\big)}\Big].
		\end{eqnarray*}
		Combining the above two inequalities, we have 
		\begin{eqnarray*}
			&\|R(x,y)-R(x',y')\| \leq  \|(x_1,x_2,y)-(x'_1,x'_2,y')\|\\
           & \times \Bigg(\E\Bigg[{\big({\sigma^{0,2}}^2+{\gamma}^2\big)\big(\big(\lambda_{21}\frac{\partial\pi^{2,*}}{\partial u}\big)^2+\big(\lambda_{22}\frac{\partial\pi^{2,*}}{\partial u}\big)^2+\big(\lambda_{22}\frac{\partial\pi^{2,*}}{\partial v}\big)^2\big)}\\&\ \ \ +{({\sigma^{0,1})}^2\big(\lambda_{11}^2\big(\frac{\partial\pi^{1,*}}{\partial u}\big)^2+\lambda_{12}^2\big(\frac{\partial\pi^{1,*}}{\partial u}\big)^2\big)}\Bigg]\Bigg)^\frac{1}{2}.
		\end{eqnarray*}
		Let $$\epsilon:=\frac{1}{6\cdot\esssup\left\{\sqrt{1+(\frac{\sigma^{0,2}}{\gamma})^2}\wedge\sqrt{1+(\frac{\gamma}{\sigma^{0,2}})^2}\right\}}.$$
		Then, for $\|\lambda\|_{\infty}<\epsilon$, we have $$\|R(x,y)-R(x',y')\|\leq L\|(x,y)-(x',y')\|$$and
        \begin{eqnarray*}
			L &&= \mathbb{E} \Big[ 			({(\lambda^{i}_{2,1})}^2+\lambda_{22}^2)\big(1+\big(\frac{\gamma}{\sigma^{0,2}}\big)^2\big)
			\frac{(\sigma^{0,2})^2}{(\sigma^{0,2})^2+\sigma^2+\gamma^2\nu e^{-\delta(\pi\gamma-\theta y)}} \\
			&&\qquad +\; \lambda_{22}^2\big(1+\big(\frac{\sigma^{0,2}}{\gamma}\big)^2\big)
			\frac{\gamma^2\nu e^{-\delta(\pi\gamma-\theta y)}}{(\sigma^{0,2})^2+\sigma^2+\gamma^2\nu e^{-\delta(\pi\gamma-\theta y)}} + \lambda_{11}^2 + \lambda_{12}^2
			\Big]\\ &&
            < \frac{5}{6} < 1,
		\end{eqnarray*}        which gives the desired result.
	\end{proof}
	Now, we give the existence of MFE in the next result.
	\begin{theorem}
		Let $(x^*_{1},x^*_{2},y^*)$ be the fixed point of Equation \eqref{eq: two player eq}.
		Then
        $(\pi^{1,*},\pi^{2,*})=\{(\pi^{1,*}_t,\pi^{2,*}_t)=(\pi^{1,*}(\lambda_{11} x^*_{1}+\lambda_{12}x^*_{2},\lambda_{12}y^*),\pi^{2,*}(\lambda_{11} x^*_{1}+\lambda_{12}x^*_{2},\lambda_{12}y^*))\}_{0\leq t\leq T}$ is the equilibrium strategy of Mean Field Game.
		Also, this equilibrium is unique in the set of all the deterministic strategy.
	\end{theorem}
	\begin{proof}
		Using the verification theorem \ref{thm: MFG optimal pi}, we know that the $(\pi^{1,*},\pi^{2,*})$ is the best response control in Problem \ref{prob: invest} satisfying the consistency condition
		\begin{eqnarray*}
			&&H^1=\lambda_{1,1}\E\[X_T^{1,\pi^{1,*}(x^*_{1})}|\F^0_T\]+\lambda_{1,2}\E\[X_T^{2,\pi^{2,*}(x^*_{2},y^*)}|\F^0_T\],\\   
			&& H^2=\lambda_{2,1}\E\[X_T^{1,\pi^{1,*}(x^*_{1})}|\F^0_T\]+\lambda_{2,2}\E\[X_T^{2,\pi^{2,*}(x^*_{2},y^*)}|\F^0_T\].
		\end{eqnarray*}
		Thus $(\pi^{1,*}_t,\pi^{2,*}_t)$ is indeed an equilibrium. For the uniqueness. We see that for every deterministic strategy $\alpha$,  $(\widebar{\alpha^1_t\sigma^{0,1}},\widebar{\alpha^2_t\sigma^{0,2}}, \widebar{\alpha^2_t\gamma})$ must be the  unique solution of Equation \eqref{eq: two player eq}. It then follows that $\alpha= (\pi^{1,*},\pi^{2,*})$.
	\end{proof}
	
	\section{Finite Population Game with $N_1+N_2$ Players}\label{sec: N1N2player}
	In this section, we study the finite player game with $N_1$ players in the population 1 and $N_2$ players in the population 2. Each agent represents a sample from the population type vector $\zeta^1$ and $\zeta^2$, that is, $N_1$ of them are sampled from $\zeta^1$ and $N_2$ of them are sampled from $\zeta^2$. We use $p=1,2$ to represent the population 1 and  the population 2, and use $i$ to indicate the $i$-th player in the population. $\zeta^{p,i}=(x_0^{p,i},\delta^{p,i},\lambda^{p,i},\mu^{p,i},\sigma^{p,i},\sigma^{0,i},\gamma^{p,i})$ is sampled from corresponding $\zeta^i$, $i=1,2$.
	
	Their wealth dynamics are given respectively by
	\begin{equation*}
		\begin{aligned}
			\ud X^{1,i}_t&=\pi^{1,i}_t\mu^{1,i}\ud t+\pi^{1,i}_t\sigma^{1,i}\ud W^{1,i}_t+\pi^{1,i}_t\sigma^{0,1,i}\ud W_t^0+\pi^{1,i}_t\gamma^{1,i}\ud \tilde{N}^{1,i}_t,\ \ 			X^{1,i}_0=x^{1,i}_0,\\
			\ud X^{2,i}_t&=\pi^{2,i}_t\mu^{2,i}\ud t+\pi^{2,i}_t\sigma^{2,i}\ud W^{2,i}_t+\pi^{2,i}_t\sigma^{0,2,i}\ud W_t^0+\pi^{2,i}_t\gamma^{2,i}\ud \tilde{N}_t, \  \  			X^{2,i}_0=x^{2,i}_0.			
		\end{aligned}		
	\end{equation*}
Here	$W^{1,i}$ and $W^{2,i}$ are idiosyncratic risk, while $W^0$ and $\tilde{N}$ are common noise.
 Suppose that $\F^0$ is the filtration generated by the common noise $W^0$ and $\tilde{N}$. 
	Again, we consider reference point process as weighted average from two populations of investors: $$H^{p,i}=\lambda^{i}_{p,1}\frac{\sum_{i=1}^{N_1}X^{1,i}_T}{N_1}+\lambda^{i}_{p,2}\frac{\sum_{i=1}^{N_2}X^{2,i}_T}{N_2},\ (p,i)\in I, $$
	where $I=\{(1,1),(1,2),\cdots,(1,N_1),(2,1),(2,2),\cdots,(2,N_2)\}$ denotes the index of $N_1+N_2$ players. 
    
	For the player $(p,i)$, the auxiliary control problem is given below.
	\begin{problem}\label{prob: invest N}
		$$\max_{\pi^{p,i}}\E[-\frac{1}{\delta^{p,i}}e^{-\delta^{p,i}(X^{p,i}_T-H^{p,i})}]$$
	for $H^{p,i}\in \F^0$.    
    \end{problem}
Our goal is to find a Nash equilibrium, defined as below.

	\begin{definition}
		We say $\{\pi^{p,i}\}_{(p,i)\in I }$ is a Nash equilibrium in the finite player game if and only if $\pi^{p,i}$ is the optimal strategy of Problem \ref{prob: invest N} with respect to $$H^{p,i}=\lambda^{i}_{p,1}\frac{\sum_{i=1}^{N_1}X^{1,i,\pi^{1,i}}_T}{N_1}+\lambda^{i}_{p,2}\frac{\sum_{i=1}^{N_2}X^{2,i,\pi^{2,i}}_T}{N_2}$$
		for every $(p,i)\in I.$
	\end{definition}
	\subsection{Solution of Problem \ref{prob: invest N}}
	First, we fix the strategy of investors except $(2,i)$ and consider Problem \ref{prob: invest N} for the investor $(2,i)$. The strategies $\{\pi^{1,j}\}_{j=1,\cdots,N}$ and $\{\pi^{2,j}\}_{j\neq i}$ are exogenous.
	We define the auxiliary process $Y^{2,i}$ by
	\begin{equation*}
		\begin{aligned}
			\ud Y^{2,i}_t=&\Big[\lambda^{i}_{2,1}\widebar{\pi\mu^1}+\lambda^{i}_{2,2}\widehat{\pi\mu}^{2,i}\Big]\ud t+\frac{\lambda^{i}_{2,1}}{N_1}\Sigma_{j=1}^{N_1}\pi^{1,j}_t\sigma^{1,j}\ud W^{1,j}_t+\frac{\lambda^{i}_{2,2}}{N_2}\Sigma_{j\neq i}\pi^{2,j}_t\sigma^{2,j}\ud W^{2,j}_t\\&+\Big[\lambda^{i}_{2,1}\widebar{\pi\sigma^{0,1}}+\lambda^{i}_{2,2}\widehat{\pi\sigma^{0,2}}^{2,i}\Big]\ud W_t^0+\lambda^{i}_{2,2}\widehat{\pi\gamma}^{2,i}\ud \tilde{N}_t+\frac{\lambda^{i}_{2,1}}{N_1}\Sigma_{j=1}^{N_1}\pi^{1,j}_t\gamma^{1,j}\ud \tilde{N}^{1,j}_t,\\
			Y^{2,i}_0=&\lambda^{i}_{2,1}\frac{\Sigma_{j=1}^{N_1}x^{1,j}_0}{N_1}+\lambda^{i}_{2,2}\frac{\Sigma_{j\neq i}x^{2,j}_0}{N_2},
		\end{aligned}		
	\end{equation*}
	where $  \widehat{\pi\mu}^{2,i}=\frac{1}{N_2}\Sigma_{j\neq i}\pi^{2,i}_t\mu^{2,i},\quad \widehat{\pi\sigma^0}^{2,i}_t=\frac{1}{N_2}\Sigma_{j\neq i}\pi^{2,i}_t\sigma^{0,2,i},\quad \widehat{\pi\gamma}^{2,i}=\frac{1}{N_2}\Sigma_{j\neq i}\pi^{2,i}_t\gamma^{2,i}$   and  $ \widebar{\pi\mu^1}=\frac{1}{N_1}\Sigma_{j=1}^{N_1}\pi^{1,i}_t\mu^{1,i},\quad \widebar{\pi\sigma^{0,1}}=\frac{1}{N_1}\Sigma_{j=1}^{N_1}\pi^{1,i}_t\sigma^{1,i}.$ 
Then, $Z_T=X^{2,i}_T-H^{2,i}=(1-\frac{\lambda^{i}_{2,2}}{N_2})X^{2,2}_T-Y^{2,i}_T$, and the auxiliary control problem \ref{prob: invest N} for the investor $(2,i)$ is formulated as follows:
	\begin{problem} \label{prob: invest N special} 
		\begin{equation*}
			\max_{\pi^{p,i}}\E\Big[-\frac{1}{\delta^{p,i}}e^{-\delta^{p,i}(X^{p,i}_T-H^{p,i})}\Big].
		\end{equation*}
	\end{problem}
	Like before, we  use the verification theorem \ref{thm: verification} to find the optimal strategy and value function. Assume $V(z,t)=-\frac{1}{\delta}e^{-\delta z}f(t)$,
	then we have \begin{equation}
		\begin{aligned}
        \min_{\pi^{2,i}_t}\big\{&-\delta^{2,i}\big[(1-\frac{\lambda^{i}_{2,2}}{N_2})\pi^{2,i}_t\mu^{2,i}-(\lambda^{i}_{2,1}\widebar{\pi\mu^1}+\lambda^{i}_{2,2}\widehat{\pi\mu}^{2,i})-((1-\frac{\lambda^{i}_{2,2}}{N_2})\pi^{2,i}_t\gamma^{2,i}-\lambda^{i}_{2,2}\widehat{\pi\gamma}^{2,i})\nu\big]\\&+\frac{1}{2}{\delta^{2,i}}^2\big[(1-\frac{\lambda^{i}_{2,2}}{N_2}){\pi^{2,i}_t}^2{\sigma^i}^2+((1-\frac{\lambda^{i}_{2,2}}{N_2})\pi^{2,i}_t\sigma^{0,2,i}-(\lambda^{i}_{2,1}\widebar{\pi\sigma^{0,1}}+\lambda^{i}_{2,2}\widehat{\pi\sigma^{0,2}}^{2,i}))^2\\&+(\frac{\lambda^{i}_{2,1}}{N_1})^2\Sigma_{j=1}^{N}(\pi^{1,j}_t\sigma^{1,j})^2+(\frac{\lambda^{i}_{2,2}}{N_2})^2\Sigma_{j\neq i}(\pi^{2,j}_t\sigma^{2,j})^2\big]\\&+\big[e^{-\delta^{2,i}((1-\frac{\lambda^{i}_{2,2}}{N_2})\pi^{2,i}_t\gamma^{2,i}-\lambda^{i}_{2,2}\widehat{\pi\gamma}^{2,i})}-1\big]\nu\big\}f(t)+f'(t)=0.
		\end{aligned}
	\end{equation}	

    \begin{theorem}\label{thm: N player invest solution}
		Let $\pi^{2,i,*}_{N_2,t}$ be the unique solution of the equation:
		\begin{equation}\label{eq: N player optimal pi 1}
			\begin{aligned}
				0=&-{\delta^{2,i}} \gamma^{2,i}\nu [e^{-{\delta^{2,i}}((1-\frac{\lambda^{i}_{2,2}}{N_2})\pi^{2,i}_t\gamma^{2,i}-\lambda^{i}_{2,2}\widehat{\pi\gamma}^{2,i})}-1] +{\delta^{2,i}}^2(1-\frac{\lambda^{i}_{2,2}}{N_2})({\sigma^{2,i}}^2+{\sigma^{0,2,i}}^2)\pi^{2,i}_t\\&-{\delta^{2,i}}^2{\sigma^{0,2,i}}(\lambda^{i}_{2,1}\widebar{\pi\sigma^{0,1}}+\lambda^{i}_{2,2}\widehat{\pi\sigma^{0,2}}^{2,i})-{\delta^{2,i}}\mu^{2,i},
			\end{aligned}
		\end{equation}
		and $V(x,t)=-\frac{1}{\delta}e^{-\delta x}f(t)$, where $f$ satisfies
		the following ODE :
	\begin{equation*}
			\begin{aligned}
				\Bigg\{&-\delta^{2,i}\Big[(1-\frac{\lambda^{i}_{2,2}}{N_2})\pi^{2,i,*}_{N_2,t}\mu^{2,i}-(\lambda^{i}_{2,1}\widebar{\pi\mu^1}+\lambda^{i}_{2,2}\widehat{\pi\mu}^{2,i})-((1-\frac{\lambda^{i}_{2,2}}{N_2})\pi^{2,i,*}_{N_2,t}\gamma^{2,i}-\lambda^{i}_{2,2}\widehat{\pi\gamma}^{2,i})\nu\Big]\\&+\frac{1}{2}{\delta^{2,i}}^2\Big[(1-\frac{\lambda^{i}_{2,2}}{N_2}){\pi^{2,i,*}_{N_2,t}}^2{\sigma^i}^2+((1-\frac{\lambda^{i}_{2,2}}{N_2})\pi^{2,i,*}_{N_2,t}\sigma^{0,2,i}-(\lambda^{i}_{2,1}\widebar{\pi\sigma^{0,1}}+\lambda^{i}_{2,2}\widehat{\pi\sigma^{0,2}}^{2,i}))^2\\&+(\frac{\lambda^{i}_{2,1}}{N_1})^2\Sigma_{j=1}^{N}(\pi^{1,j}_t\sigma^{1,j})^2+(\frac{\lambda^{i}_{2,2}}{N_2})^2\Sigma_{j\neq i}(\pi^{2,j}_t\sigma^{2,j})^2\Big]\\&+\Big[e^{-\delta^{2,i}((1-\frac{\lambda^{i}_{2,2}}{N_2})\pi^{2,i,*}_{N_2,t}\gamma^{2,i}-\lambda^{i}_{2,2}\widehat{\pi\gamma}^{2,i})}-1\Big]\nu\Bigg\}f(t)+f'(t)=0.
			\end{aligned}
		\end{equation*}
		Then $\pi^{2,i,*}_{N_2,t}$ and $V(x,t)$ are the optimal strategy and value function of Problem \ref{prob: invest N special}.
	\end{theorem}
	\begin{proof} Because
		the proof is similar to that of Theorem \ref{thm: MFG optimal pi} , we omit it here.
	\end{proof}
	
	It is worth mentioning that Equation \eqref{eq: N player optimal pi 1} is equivalent to the following equation:
	$$\begin{aligned}
		0=&-{\delta^{2,i}} \gamma^{2,i}\nu \Big[e^{-{\delta^{2,i}}\(\pi^{2,i}\gamma^{2,i}-\lambda^{i}_{2,2} \widebar{\pi^2\gamma}\)}-1\Big]+{\delta^{2,i}}^2\Big((1-\frac{\lambda^{i}_{2,2}}{N_2}){\sigma^{2,i}}^2+{\sigma^{0,2,i}}^2\Big)\pi^{2,i}\\&-{\delta^{2,i}}^2{\sigma^{0,2,i}}(\lambda^{i}_{2,1}\widebar{\pi^1\sigma^{0,1}}_t+\lambda^{i}_{2,2}\widebar{\pi^2\sigma^{0,2}}_t)-{\delta^{2,i}}\mu^{2,i},
	\end{aligned}$$
	where the definition of $\widebar{\cdot}$ is same as the previous notation.
	This conversion is useful for us to find the equilibrium as we change the asymmetric mean $\widehat{\pi\sigma}^{2,i}_t$ and $\widehat{\pi\gamma}^{2,i}$ into symmetric mean $\widebar{\pi^2\sigma^{0,2}}$ and $\widebar{\pi^2\gamma}$.
	From now on, we define $g^{2,i}$ by
	$$\begin{aligned}
		g^{2,i}(\pi,u,v)=&-{\delta^{2,i}} \gamma^{2,i}\nu \[e^{-{\delta^{2,i}}\(\pi\gamma^{2,i}-v\)}-1\]+{\delta^{2,i}}^2\Big((1-\frac{\lambda^{i}_{2,2}}{N_2}){\sigma^{2,i}}^2+{\sigma^{0,2,i}}^2\Big)\pi\\&-{\delta^{2,i}}^2{\sigma^{0,2,i}}u-{\delta^{2,i}}\mu^{2,i}.
	\end{aligned}$$
	Then we have the next result.
	\begin{lemma}\label{lem: 2N equiv eq}
		$\pi^{2,i,*}_{N_2}$ defined in Theorem \ref{thm: N player invest solution} is the unique solution of the following equation with respect to $\pi^{2,i}$:
		$$\begin{aligned}
			&g^{2,i}\(\pi^{2,i},\lambda^{i}_{2,1}\widebar{\pi^1\sigma^{0,1}}+\lambda^{i}_{2,2}\widebar{\pi^2\sigma^{0,2}},\lambda^{i}_{2,2}\widebar{\pi^2\gamma}\)\\&=g^{2,i}\Big(\pi^{2,i},\lambda^{i}_{2,1}\widebar{\pi^1\sigma^{0,1}}+\lambda^{i}_{2,2}\widehat{\pi\sigma}^{2,i}_t+\frac{\lambda^{i}_{2,2}}{N_2}\pi^{2,i}\sigma^{0,2,i},\lambda^{i}_{2,2}\widehat{\pi\gamma}^{2,i}+\frac{\lambda^{i}_{2,2}}{N_2}\pi^{2,i}\gamma^{2,i}\Big)=0.
		\end{aligned}$$
	\end{lemma}
	Now we  use the same procedure to obtain the optimal strategy of the investor $(1,i)$.
	\begin{theorem}
		The optimal strategy $\pi^{1,i,*}_{N_1}$ for the investor $(1,i)$ solving Problem \ref{prob: invest N} is the unique solution of the following equation:
		$$
		\begin{aligned}
			0=&-{\delta^{1,i}} \gamma^{1,i}\nu \Big[e^{-{\delta^{1,i}}((1-\frac{\lambda^{i}_{1,1}}{N_1})\pi^{1,i}_t\gamma^i_t)}-1\Big]+{\delta^{1,i}}^2(1-\frac{\lambda^{i}_{1,1}}{N_1})\({\sigma^{1,i}}^2+{\sigma^{0,1,i}}^2\)\pi^{1,i}_t\\&-{\delta^{1,i}}^2{\sigma^{0,2,i}}\Big(\lambda^{i}_{1,2}\widebar{\pi\sigma^{0,2}}_t+\lambda^{i}_{1,1}\widehat{\pi\sigma^{0,1}}^{1,i}_t\Big)-{\delta^{1,i}}\mu^{1,i}.
		\end{aligned}
		$$
		Define $$\begin{aligned}
			g^{1,i}(\pi,u,v)=&-{\delta^{1,i}} \gamma^{1,i}\nu \Big[e^{-{\delta^{1,i}}((1-\frac{\lambda^{i}_{1,1}}{N_1})\pi\gamma^{2,i}-v)}-1\Big]+{\delta^{1,i}}^2\Big((1-\frac{\lambda^{i}_{1,1}}{N_1}){\sigma^{1,i}}^2+{\sigma^{0,1,i}}^2\Big)\pi\\&-{\delta^{1,i}}^2{\sigma^{0,1,i}}u-{\delta^{1,i}}\mu^{1,i}.
		\end{aligned}$$
		Then $\pi^{1,i,*}_{N_1}$ is also the unique solution of the following equation with respect to $\pi^{1,i}$:
		$$\begin{aligned}
			&g^{1,i}\(\pi^{1,i},\lambda^{i}_{1,1}\widebar{\pi^1\sigma^{0,1}}+\lambda^{i}_{1,2}\widebar{\pi^2\sigma^{0,2}},0\)\\&=g^{1,i}\Big(\pi^{1,i},\lambda^{i}_{1,1}\widehat{\pi\sigma}^{1,i}_t+\frac{\lambda^{i}_{1,1}}{N_1}\pi^{1,i}\sigma^{0,1,i}+\lambda^{i}_{1,2}\widebar{\pi^2\sigma^{0,2}},0\Big)=0.
		\end{aligned}$$
	\end{theorem}
	\subsection{Finding Nash equilibrium} 
	The goal of this subsection is to characterize the Nash equilibrium of the $N_1+N_2$-player game by employing the fixed point argument.
	
	Denote $x_1=\widebar{\pi^1\sigma^{0,1}},x_2=\widebar{\pi^2\sigma^{0,2}},y=\widebar{\pi^2\gamma}$, and  $\mathbf{N}=(N_1,N_2)$  the vector numbers of the two populations.
	Then we have the following theorem on the Nash equilibrium. 
	\begin{theorem}\label{thm: $N_1+N_2$-player equi verification}
		Let $\pi^{1,i,*}_{N_1}(u,0)$ and $\pi^{2,i,*}_{N_2}(u,v)$ be the unique solution of equations\\ $g^{1,i}(\pi,u,0)=0$ and $g^{2,i}(\pi,u,v)=0$, respectively.
		If the  equation system
		\begin{equation}\label{eq: equi of 2*N}
			\left\{\begin{aligned}
				&\Sigma_{i =1}^{N_1}\frac{\pi^{1,i,*}_{N_1}(\lambda^{i}_{1,1} x_1+\lambda^{i}_{1,2}x_2,0)\sigma^{0,1,i}}{N_1}=x_1,\\
				&\Sigma_{i =1}^{N_2}\frac{\pi^{2,i,*}_{N_2}(\lambda^{i}_{2,1} x_1+\lambda^{i}_{2,2}x_2,\lambda^{i}_{2,2}y)\sigma^{0,2,i}}{N_2}=x_2,\\
				&\Sigma_{i =1}^{N_2}\frac{\pi^{2,i,*}_{N_2}(\lambda^{i}_{2,1} x_1+\lambda^{i}_{2,2}x_2,\lambda^{i}_{2,2}y)\gamma_i}{N_2}=y.\\
			\end{aligned}\right.
		\end{equation}
		has a fixed point $(x^*_{1,\mathbf{N}},x^*_{2,\mathbf{N}},y^*_{\mathbf{N}})$, then
		\begin{equation}\label{eq: strategy 2N}
			\{(\pi^{1,i,*}_{N_1},\pi^{2,i,*}_{N_2})=\(\pi^{1,i,*}_{N_1}(\lambda^{i}_{1,1} x^*_{1,\mathbf{N}}+\lambda^{i}_{1,2}x^*_{2,\mathbf{N}},0),\pi^{2,i,*}_{N_2}(\lambda^{i}_{2,1} x^*_{1,\mathbf{N}}+\lambda^{i}_{2,2}x^*_{2,\mathbf{N}},\lambda^{i}_{2,2}y^*_{\mathbf{N}})\)\}_{1\leq i\leq N}
		\end{equation} is a Nash equilibrium in the $N_1+N_2$ player game.
	\end{theorem}
	\begin{proof}
		Suppose that  $(x^*_{1,\mathbf{N}},x^*_{2,\mathbf{N}},y^*_{\mathbf{N}})$ is a solution of Equation (\ref{eq: equi of 2*N}). Then we only need to prove that the strategy \eqref{eq: strategy 2N} is a Nash equilibrium strategy. First, we suppose that all the investors except $(2,i)$ choose this strategy. 
		
		According to Theorem \ref{thm: N player invest solution} and Lemma \ref{lem: 2N equiv eq}, we know that the optimal strategy of the  investor $(2,i)$ is the unique solution of the following equation
		\begin{equation}\label{eq: g2i1}
			g^{2,i}\Big(\pi^{2,i},\lambda^{i}_{2,1}\widebar{\pi^{1,*}\sigma^{0,1}}+\lambda^{i}_{2,2}\widehat{\pi^{*}\sigma}^{2,i}+\frac{\lambda^{i}_{2,2}}{N_2}\pi^{2,i}\sigma^{0,2,i},\lambda^{i}_{2,2}\widehat{\pi^{*}\gamma}^{2,i}+\frac{\lambda^{i}_{2,2}}{N_2}\pi^{2,i}\gamma^{2,i}\Big)=0.
		\end{equation}
		Meanwhile, based on the definition of $(\pi^{1,i,*}_{N_1},\pi^{2,i,*}_{N_2})$, we have $$g^{2,i}(\pi^{2,i,*}_{N_2}(\lambda^{i}_{2,1} x^*_{1,\mathbf{N}}+\lambda^{i}_{2,2}x^*_{2,\mathbf{N}},\lambda^{i}_{2,2}y^*_{\mathbf{N}}),\lambda^{i}_{2,1} x^*_{1,\mathbf{N}}+\lambda^{i}_{2,2}x^*_{2,\mathbf{N}},\lambda^{i}_{2,2}y^*_{\mathbf{N}})=0.$$
		Thanks to Equation \eqref{eq: equi of 2*N}, the last equation  becomes $$g^{2,i}\(\pi^{2,i,*}_{N_2},\lambda^{i}_{2,1}\widebar{\pi^{1,*}\sigma^{0,1}}+\lambda^{i}_{2,2}\widebar{\pi^{2,*}\sigma^{0,2}},\lambda^{i}_{2,2}\widebar{\pi^{2,*}\gamma}\)=0.$$
		Thus 
		\begin{equation}\label{eq: g2i2}
			g^{2,i}\Big(\pi^{2,i,*}_{N_2},\lambda^{i}_{2,1}\widebar{\pi^{1,*}\sigma^{0,1}}+\lambda^{i}_{2,2}\widehat{\pi^{*}\sigma}^{2,i}+\frac{\lambda^{i}_{2,2}}{N_2}\pi^{2,i,*}_{N_2}\sigma^{0,2,i},\lambda^{i}_{2,2}\widehat{\pi^{*}\gamma}^{2,i}+\frac{\lambda^{i}_{2,2}}{N_2}\pi^{2,i,*}_{N_2}\gamma^{2,i}\Big)=0.
		\end{equation} 
		Comparing Equations \eqref{eq: g2i1} and \eqref{eq: g2i2}, we know that the optimal strategy of  the investor $(2,i)$ must be $\pi^{2,i,*}_{N_2}$.
		Similarly, $\pi^{1,i,*}_{N_1}$ is the optimal strategy of the investor $(1,i)$ when all the other investors choose strategy \eqref{eq: strategy 2N}. Therefore Strategy \eqref{eq: strategy 2N} is the Nash equilibrium strategy.
	\end{proof}
	Now, the remaining task is to establish the existence of unique solution to equation system (\ref{eq: equi of 2*N}). We can adopt the same approach for MFG in Subsection \ref{MFG}. 

	Based on  Implicit function theorem, 
    $$\frac{\partial\pi^{1,i,*}_{N_1}(u,v)}{\partial x}=-\frac{\frac{\partial g^{1,i}}{\partial u}}{\frac{\partial g_i}{\partial \pi}}=\frac{\sigma^{0,1,i}}{({\sigma^{0,1,i}}^2+(1-\frac{\lambda^{i}_{1,1}}{N_1}){\sigma^{1,i}}^2)+{\gamma^{1,i}}^2\nu e^{-\delta(1-\frac{\lambda^{i}_{1,1}}{N_1})\(\pi\gamma^{1,i}\)}},$$
	$$\frac{\partial\pi^{1,i,*}_{N_1}(u,v)}{\partial v}=-\frac{\frac{\partial g^{2,i}}{\partial v}}{\frac{\partial g_i}{\partial \pi}}=\frac{\gamma^{1,i}\nu e^{-\delta(1-\frac{\lambda^{i}_{1,1}}{N_1})\(\pi\gamma^{1,i}\)}}{({\sigma^{0,1,i}}^2+(1-\frac{\lambda^{i}_{1,1}}{N_1}){\sigma^{1,i}}^2)+{\gamma^{1,i}}^2\nu e^{-\delta(1-\frac{\lambda^{i}_{1,1}}{N_1})\(\pi\gamma^{1,i}\)}}.$$
	$$\frac{\partial\pi^{2,i,*}_{N_2}(u,v)}{\partial x}=-\frac{\frac{\partial g^{2,i}}{\partial u}}{\frac{\partial g_i}{\partial \pi}}=\frac{\sigma^{0,2,i}}{({\sigma^{0,2,i}}^2+(1-\frac{\lambda^{i}_{2,2}}{N_2}){\sigma^{2,i}}^2)+{\gamma^{2,i}}^2\nu e^{-\delta\(\pi\gamma^{2,i}-v\)}},$$
	$$\frac{\partial\pi^{2,i,*}_{N_2}(u,v)}{\partial v}=-\frac{\frac{\partial g^{2,i}}{\partial v}}{\frac{\partial g_i}{\partial \pi}}=\frac{\gamma^{2,i}\nu e^{-\delta\(\pi\gamma^{2,i}-v\)}}{({\sigma^{0,2,i}}^2+(1-\frac{\lambda^{i}_{2,2}}{N_2}){\sigma^{2,i}}^2)+{\gamma^{2,i}}^2\nu e^{-\delta\(\pi\gamma^{2,i}-v\)}}.$$
	For small $\theta$, we have a unique equilibrium:
	\begin{theorem}
		There exists an  $\epsilon>0$ depending on parameters $\sigma^i,\sigma^{i,0},\gamma^i$ such that if $\|\lambda\|_{\infty}<\epsilon$, equation system \eqref{eq: equi of 2*N} admits a unique solution $(x_1,x_2,y)\in \R^2$. In fact, we could use same $\epsilon$ in Theorem \ref{thm: MFG existent of equi}.
	\end{theorem}
	\begin{proof}
		Define
		$$
		\begin{aligned}
			R_{1,N_1}(x_1,x_2,y)=&\Sigma_{i =1}^{N_1}\frac{\pi^{1,i,*}_{N_1}(\lambda^{i}_{1,1} x_1+\lambda^{i}_{1,2}x_2,0)\sigma^{0,1,i}}{N_1},\\
			R_{2,N_2}(x_1,x_2,y)=&\Sigma_{i =1}^{N_2}\frac{\pi^{2,i,*}_{N_2}(\lambda^{i}_{2,1} x_1+\lambda^{i}_{2,2}x_2,\lambda^{i}_{2,2}y)\sigma^{0,2,i}}{N_2},\\
			R_{3,N_2}(x_1,x_2,y)=&\Sigma_{i =1}^{N_2}\frac{\pi^{2,i,*}_{N_2}(\lambda^{i}_{2,1} x_1+\lambda^{i}_{2,2}x_2,\lambda^{i}_{2,2}y)\gamma_i}{N_2}
		\end{aligned}
		$$and 
		$$R_{\mathbf{N}}(x_1,x_2,y)=\(\begin{matrix}
			R_{1,N_1}(x_1,x_2,y),\\
			R_{2,N_2}(x_1,x_2,y),\\
			R_{3,N_2}(x_1,x_2,y).\\
		\end{matrix}\).$$
Similar to Theorem \ref{thm: MFG existent of equi}, we only need to prove that $R_{\mathbf{N}}$ is a contraction mapping from $\R^3$ to $\R^3$.
		For every $(x_1,x_2,y), (x'_1,x'_2,y')\in \R^3$, 
		$$\begin{aligned}
			&R_{2,N_2}(x_1,x_2,y)-R_{2,N_2}(x'_1,x'_2,y')\\
			=&\Sigma_{i =1}^{N_2}\frac{1}{N_2}\Big[\sigma^{0,2,i}\[{\lambda^{i}_{2,1}}(x_1-x'_1)+{\lambda^{i}_{2,2}}(x_2-x'_2)\]\frac{\partial\pi^{2,i,*}_{N_2}}{\partial u}+\sigma^{0,2,i}\[{\lambda^{i}_{2,2}}(y-y')\]\frac{\partial\pi^{2,i,*}_{N_2}}{\partial v}\Big]\\
			\leq &\|(x_1,x_2,y)-(x'_1,x'_2,y')\|\Sigma_{i =1}^{N_2}\frac{1}{N_2}\Big[\sqrt{({\sigma^{0,2,i}})^2\Pi^{2,i}}\Big],
		\end{aligned}$$
        where 
        $$\Pi^{2,i}={(\lambda^{i}_{2,1})}^2\Big(\frac{\partial\pi^{2,i,*}_{N_2}}{\partial u}\Big)^2+{(\lambda^{i}_{2,2})}^2\Big(\frac{\partial\pi^{2,i,*}_{N_2}}{\partial u}\Big)^2+{(\lambda^{i}_{2,2})}^2\Big(\frac{\partial\pi^{2,i,*}_{N_2}}{\partial v}\Big)^2$$
        evaluate at $(u,v)={\lambda^{i}_{2,1}}(x'_1+\theta(x_1-x'_1))+{\lambda^{i}_{2,2}}(x'_2+\theta(x_2-x'_2)),{\lambda^{i}_{2,2}}(y'+\theta(y-y'))$.
        
		Using Jensen's inequality, we have 
		{\small
			$$\begin{aligned}
				&\(R_{2,N_2}(x_1,x_2,y)-R_{2,N_2}(x'_1,x'_2,y')\)^2\\
				\leq&\|(x_1,x_2,y)-(x'_1,x'_2,y')\|^2\Big[\Sigma_{i =1}^{N_2}\frac{1}{N_2}\sqrt{({\sigma^{0,2,i}})^2\Pi^{2,i}}\Big]^2\\
				\leq&\|(x_1,x_2,y)-(x'_1,x'_2,y')\|^2\Sigma_{i =1}^{N_2}\frac{1}{N_2}\[({\sigma^{0,2,i}})^2\Pi^{2,i}\].
			\end{aligned}$$}
            Repeating the process, we get 
		{\small
			$$\begin{aligned}
				&\(R_{1,N_1}(x_1,x_2,y)-R_{1,N_1}(x'_1,x'_2,y')\)^2\\
				\leq&\|(x_1,x_2,y)-(x'_1,x'_2,y')\|^2\Sigma_{i =1}^{N_1}\frac{1}{N_1}\Big[({\sigma^{0,1,i}})^2\Big({(\lambda^{i}_{1,1})}^2\Big(\frac{\partial\pi^{1,i,*}_{N_1}}{\partial u}\Big)^2+{(\lambda^{i}_{1,2})}^2\Big(\frac{\partial\pi^{1,i,*}_{N_1}}{\partial u}\Big)^2\Big)\Big],
			\end{aligned}$$}and
		{\small
			$$\begin{aligned}
				\(R_{3,N_2}(x_1,x_2,y)-R_{3,N_2}(x'_1,x'_2,y')\)^2
				\leq\|(x_1,x_2,y)-(x'_1,x'_2,y')\|^2\Sigma_{i =1}^{N_2}\frac{1}{N_2}\Big[({\gamma^{2,i}})^2\Pi^{2,i}\Big].
			\end{aligned}$$}
		Define $$\epsilon=\frac{1}{6\cdot\esssup\left\{\sqrt{1+(\frac{\sigma^{0,2}}{\gamma})^2}\wedge\sqrt{1+(\frac{\gamma}{\sigma^{0,2}})^2}\right\}}.$$
		Then, for $\|\lambda\|_{\infty}<\epsilon$, we have $$\|R_{\mathbf{N}}(x_1,x_2,y)-R_{\mathbf{N}}(x'_1,x'_2,y')\|\leq L\|(x_1,x_2,y)-(x'_1,x'_2,y')\|, $$ where 
        $$\begin{aligned}
			&L = \\
            &\Sigma_{i =1}^{N_2}\frac{1}{N_2}\Big[ 
			({(\lambda^{i}_{2,1})}^2+{(\lambda^{i}_{2,2})}^2)\Big(1+(\frac{\gamma^{2,i}}{\sigma^{0,2}})^2\Big)
			\frac{(\sigma^{0,2,i})^2}{({\sigma^{0,2,i}}^2+(1-\frac{\lambda^{i}_{2,2}}{N_2}){\sigma^{2,i}}^2)+{\gamma^{2,i}}^2\nu e^{-\delta\(\pi\gamma^{2,i}-v\)}} \\
			&+{(\lambda^{i}_{2,2})}^2\Big(1+(\frac{\sigma^{0,2,i}}{\gamma^{2,i}})^2\Big)
			\frac{{\gamma^{2,i}}^2\nu e^{-\delta\(\pi\gamma^{2,i}-v\)}}{({\sigma^{0,2,i}}^2+(1-\frac{\lambda^{i}_{2,2}}{N_2}){\sigma^{2,i}}^2)+{\gamma^{2,i}}^2\nu e^{-\delta\(\pi\gamma^{2,i}-v\)}}\Big] \\&+ {(\lambda^{i}_{1,1})}^2 + {(\lambda^{i}_{1,2})}^2
			< \frac{5}{6} < 1,
		\end{aligned}$$
        which yields the desired result.
	\end{proof}
    
    \section{Relation between MFG and the $N_1+N_2$-Player Game}\label{sec: convergence}
	In this section, we discuss the relationship between the MFE for the MFG and the Nash equilibrium for the $N_1+N_2$-player game. We show that the MFE is exactly the limit of the Nash equilibrium as $N_1,N_2\rightarrow\infty$.  
	
	Denote $x^*_{1,\mathbf{N}}$, $x^*_{2,\mathbf{N}}$ and $y^*_{\mathbf{N}}$ as the the equilibrium mean of the N player game, which is the solution of Equation \eqref{eq: two player eq}. Recall that the optimal solution of $N_1+N_2$-player game is $\(\pi^{1,i,*}_{N_1}(\lambda^{i}_{1,1} x^*_{1,\mathbf{N}}+\lambda^{i}_{1,2}x^*_{2,\mathbf{N}},0),\pi^{2,i,*}_{N_2}(\lambda^{i}_{2,1} x^*_{1,\mathbf{N}}+\lambda^{i}_{2,2}x^*_{2,\mathbf{N}},\lambda^{i}_{2,2}y^*_{\mathbf{N}})\)$, where the $\pi^{1,i,*}_{N_1}(u,v)$ and $\pi^{2,i,*}_{N_2}(u,v)$ are the unique solutions of the following equations:
	$$\begin{aligned}
		g^{1,i}(\pi,u,v)=&-{\delta^{1,i}} \gamma^{1,i}\nu \Big[e^{-{\delta^{1,i}}\big((1-\frac{\lambda^{i}_{1,1}}{N_1})\pi\gamma^{2,i}-v\big)}-1\Big]+{\delta^{1,i}}^2\Big((1-\frac{\lambda^{i}_{1,1}}{N_1}){\sigma^{1,i}}^2+{\sigma^{0,1,i}}^2\Big)\pi\\&-{\delta^{1,i}}^2{\sigma^{0,1,i}}u-{\delta^{1,i}}\mu^{1,i}=0,\\
		g^{2,i}(\pi,u,v)=&-{\delta^{2,i}} \gamma^{2,i}\nu \Big[e^{-{\delta^{2,i}}\big(\pi\gamma^{2,i}-v\big)}-1\Big]+{\delta^{2,i}}^2\Big((1-\frac{\lambda^{i}_{2,2}}{N_2}){\sigma^{2,i}}^2+{\sigma^{0,2,i}}^2\Big)\pi\\&-{\delta^{2,i}}^2{\sigma^{0,2,i}}u-{\delta^{2,i}}\mu^{2,i}=0.
	\end{aligned}$$
	The functions $g^{1,i}(\pi,u,v)$ and $g^{2,i}(\pi,u,v)$ vary in $i$ only through the type vector $\zeta^{1,i}$ and $\zeta^{2,i}$. Thus, we can regard the $\pi^{p,i,*}_{N_p}(u,v)$ as a deterministic function of $\(\zeta^{p,i},u,v\)$, i.e., $$\pi^{p,i,*}_{N_p}(u,v)=\pi^{p,*}_{N_p}(\zeta^{p,i},u,v),$$
	where $\pi^{p,*}_{N_p}(\zeta,u,v)$ is the solution of the following equations (assume $\zeta=(\xi,\delta,\lambda,\mu,\sigma,\sigma^0,\gamma)$):
	$$\begin{aligned}
		g^{1}_{N_1}(\pi,u,v,\zeta)=&-{\delta} \gamma\nu \Big[e^{-{\delta}(1-\frac{\lambda_{1}}{N_1})\(\pi\gamma\)}-1\Big]+{\delta}^2\Big((1-\frac{\lambda_{2}}{N_1}){\sigma}^2+{\sigma^{0}}^2\Big)\pi-{\delta^{2}}^2{\sigma^{0}_t}u-{\delta}\mu=0,\\
		g^{2}_{N_2}(\pi,u,v,\zeta)=&-{\delta} \gamma\nu \Big[e^{-{\delta}\(\pi\gamma-v\)}-1\Big]+{\delta}^2\Big((1-\frac{\lambda_{2}}{N_2}){\sigma}^2+{\sigma^{0}}^2\Big)\pi-{\delta^{2}}^2{\sigma^{0}_t}u-{\delta}\mu=0.
	\end{aligned}$$
	Then we have the following theorem.
	\begin{theorem}
		Regarding the  $\pi^{p,i,*}_{N_p}(u,v)$ and $\pi^{p,*}_{N_p}(u,v)$ as deterministic functions of\\ $\(\zeta^{p,i},u,v\)$ and $\(\zeta^p,u,v\)$:
		$$\pi^{p,i,*}_{N_p}(u,v)=\pi^{p,*}_{N_p}(\zeta^{p,i},u,v),\ \pi^{p,*}(u,v)=\pi^{p,*}(\zeta^{p},u,v),$$
	we have
		$$\pi^{p,*}_{N_p}\xrightarrow{N_p\to \infty}\pi^{p,*}\  \text{ in }\R^8.$$
		Also, fix the player $(p, i)$ and its projected type vector $\zeta^{p,i}$, then the Nash equilibrium strategy for the player $i$ satisfies $$\pi^{1,*}_{N_1}(\zeta^{1,i},\lambda^{i}_{1,1} x^*_{1,\mathbf{N}}+\lambda^{i}_{1,2}x^*_{2,\mathbf{N}},0)\xrightarrow{N_1,N_2\to \infty}\pi^{1,*}(\zeta^{1,i},\lambda^{i}_{1,1} x^*_{1}+\lambda^{i}_{1,2}x^*_{2},0),$$
		$$\pi^{2,*}_{N_2}(\zeta^{2,i},\lambda^{i}_{2,1} x^*_{1,\mathbf{N}}+\lambda^{i}_{2,2}x^*_{2,\mathbf{N}},\lambda^{i}_{2,2}y^*_{\mathbf{N}})\xrightarrow{N_1,N_2\to \infty}\pi^{2,*}(\zeta^{2,i},\lambda^{i}_{2,1} x^*_{1}+\lambda^{i}_{2,2}x^*_{2},\lambda^{i}_{2,2}y^*),$$ which indicates that the Nash equilibrium of $\mathbf{N}$ player game converges to the MFE.
	\end{theorem}
\begin{proof} First, it is easy to see that $\pi^{p,*}_{N_p}(\zeta,u,v)$ is monotonic with respect to $N$. Thus, we can get the convergence in $\R^8$ based on  the definition of $\pi^{p,*}_{N_p}$ and $\pi^{p,*}$. 
		
		To prove the convergence result, we have to first verify that $(x^*_{1,\mathbf{N}},x^*_{2,\mathbf{N}},y^*_{\mathbf{N}})\xrightarrow{N_1,N_2\to \infty}(x^*_{1},x^*_{2},y^*)$. 
		Recall that $(x^*_{1,\mathbf{N}},x^*_{2,\mathbf{N}},y^*_{\mathbf{N}})$ and $(x^*_1,x^*_2,y^*)$ are the fixed points of contracting map $R_{\mathbf{N}}$ and $R$, respectively. Then
		$$\begin{aligned}
			&(x^*_{1,\mathbf{N}},x^*_{2,\mathbf{N}},y^*_{\mathbf{N}})-(x^*_1,x^*_2,y^*)\\=&R_{\mathbf{N}}(x^*_{1,\mathbf{N}},x^*_{2,\mathbf{N}},y^*_{\mathbf{N}})-R(x^*_1,x^*_2,y^*)\\
			=&R_{\mathbf{N}}(x^*_{1,\mathbf{N}},x^*_{2,\mathbf{N}},y^*_{\mathbf{N}})-R_{\mathbf{N}}(x^*_1,x^*_2,y^*)+R_{\mathbf{N}}(x^*_1,x^*_2,y^*)-R(x^*_1,x^*_2,y^*).
		\end{aligned}$$
		Using the convergence of $\pi^{i,*}_N(\sigma^i,\sigma^{0,i},\gamma^i,u,v)$ and Law of Large Numbers, we have $$R_{\mathbf{N}}(x_1,x_2,y)\xrightarrow{N_1,N_2\to \infty}R(x_1,x_2,y)\  \ in\ \R^3.$$ Thus, for  $\epsilon>0$, there exists some $N_0$ such that for  $\forall$ \  $N_1,N_2>N_0$, we have
		$$\begin{aligned}
			\|(x^*_{1,\mathbf{N}},x^*_{2,\mathbf{N}},y^*_{\mathbf{N}})-(x^*_1,x^*_2,y^*)\|\leq& \|R_{\mathbf{N}}(x^*_{1,\mathbf{N}},x^*_{2,\mathbf{N}},y^*_{\mathbf{N}})-R_{\mathbf{N}}(x^*_1,x^*_2,y^*)\|+\epsilon\\
			\leq &L\|(x^*_{1,\mathbf{N}},x^*_{2,\mathbf{N}},y^*_{\mathbf{N}})-(x^*_1,x^*_2,y^*)\|+\epsilon
		\end{aligned}
		$$ with $L<1$,  yielding $(x^*_{1,\mathbf{N}},x^*_{2,\mathbf{N}},y^*_{\mathbf{N}})\xrightarrow{N_1,N_2\to \infty}(x^*_1,x^*_2,y^*)$. 
		Then
		$$\begin{aligned}
			&|\pi^{2,*}_{N_2}(\zeta^{2,i},\lambda^{i}_{2,1} x^*_{1,\mathbf{N}}+\lambda^{i}_{2,2}x^*_{2,\mathbf{N}},\lambda^{i}_{2,2}y^*_{\mathbf{N}})-\pi^{2,*}(\zeta^{2,i},\lambda^{i}_{2,1} x^*_{1}+\lambda^{i}_{2,2}x^*_{2},\lambda^{i}_{2,2}y^*)|\\=&|\pi^{2,*}_{N_2}(\zeta^{2,i},\lambda^{i}_{2,1} x^*_{1,\mathbf{N}}+\lambda^{i}_{2,2}x^*_{2,\mathbf{N}},\lambda^{i}_{2,2}y^*_{\mathbf{N}})-\pi^{2,*}_{N_2}(\zeta^{2,i},\lambda^{i}_{2,1} x^*_{1}+\lambda^{i}_{2,2}x^*_{2},\lambda^{i}_{2,2}y^*)\\&+\pi^{2,*}_{N_2}(\zeta^{2,i},\lambda^{i}_{2,1} x^*_{1}+\lambda^{i}_{2,2}x^*_{2},\lambda^{i}_{2,2}y^*)-\pi^{2,*}(\zeta^{2,i},\lambda^{i}_{2,1} x^*_{1}+\lambda^{i}_{2,2}x^*_{2},\lambda^{i}_{2,2}y^*)|\\
			\leq &L\|(x^*_{2,\mathbf{N}},y^*_{N})-(x^*,y^*)\|\\
			&+|\pi^{2,*}_{N_2}(\zeta^{2,i},\lambda^{i}_{2,1} x^*_{1}+\lambda^{i}_{2,2}x^*_{2},\lambda^{i}_{2,2}y^*)-\pi^{2,*}(\zeta^{2,i},\lambda^{i}_{2,1} x^*_{1}+\lambda^{i}_{2,2}x^*_{2},\lambda^{i}_{2,2}y^*)|.
		\end{aligned}$$
		Letting $N_1,N_2\to \infty$, we obtain $$\pi^{2,*}_{N_2}(\zeta^{2,i},\lambda^{i}_{2,1} x^*_{1,\mathbf{N}}+\lambda^{i}_{2,2}x^*_{2,\mathbf{N}},\lambda^{i}_{2,2}y^*_{\mathbf{N}})\xrightarrow{N_1,N_2\to \infty}\pi^{2,*}(\zeta^{2,i},\lambda^{i}_{2,1} x^*_{1}+\lambda^{i}_{2,2}x^*_{2},\lambda^{i}_{2,2}y^*).$$
		Similarly, it holds that
        $$\pi^{1,*}_{N_1}(\zeta^{1,i},\lambda^{i}_{1,1} x^*_{1,\mathbf{N}}+\lambda^{i}_{1,2}x^*_{2,\mathbf{N}},0)\xrightarrow{N_1,N_2\to \infty}\pi^{1,*}(\zeta^{1,i},\lambda^{i}_{1,1} x^*_{1}+\lambda^{i}_{1,2}x^*_{2},0).$$
	\end{proof}
    In conclusion, we have proved that from both perspectives—whether viewing the equilibrium strategy as a random variable or as a function of the market parameters and the equilibrium relative point $(u,v)$, the MFE is indeed the limit of the Nash equilibirum in the two-population finite-player game.

	\section{Discussions on MFE and the Impact of Poisson jump risk}\label{sec: numerical}
	In this section, we will first discuss some degenerate cases in one population MFG and analyze the equilibrium behavior under Poisson idiosyncratic and common noise. Then, for the general MFG with two heterogeneous populations, we present some numerical plots to illustrate some quantitative behavior and sensitivity analysis of MFE with respect to some market parameters.
	
	\subsection{MFG and N player game with one population}\label{subsec: only brownian}
	 {Lacker and Zariphopoulou \cite{LD2019mean}} considered MFG and N player game with Brownian common noise. They have only one population with type vector $\zeta=(\xi,\delta,\theta,\mu,\sigma,\sigma^{0,1})$. $\theta$ is a one dimensional random variable, representing the awareness of the relative performance. When consider $\lambda_{1,1}=\theta$, $\lambda_{1,2}=0$ and $\gamma^1=0$, our second population will degenerate and the problem of first population is exactly the same with {Lacker and Zariphopoulou \cite{LD2019mean}}. When $\lambda_{1,2}=0$, the population 1 will only consider their own mean of final wealth, ignoring the performance of the population 2.
	
	Under this setting, we note that the first equation in \eqref{eq: two player eq} is equal to 
 	$$\E\Big[\frac{{\sigma^{0,1}}^2}{{\sigma^{1}}^2+{\sigma^{0,1}}^2}(\lambda_{11}x_1+\lambda_{12}x_2)+\frac{\mu^1\sigma^{0,1}}{\delta^1({\sigma^{1}}^2+{\sigma^{0,1}}^2)}\Big]=x_1.$$
	Thus, the relationship between $x_1$ and $x_2$ becomes linear:
	$$x_1=\frac{\E\big[\frac{{\sigma^{0,1}}^2\lambda_{12}}{{\sigma^{1}}^2+{\sigma^{0,1}}^2}\big]}{1-\E\big[\frac{{\sigma^{0,1}}^2\lambda_{11}}{{\sigma^{1}}^2+{\sigma^{0,1}}^2}\big]}x_2+\frac{\E\big[\frac{\mu^1\sigma^{0,1}}{\delta^1\({\sigma^{1}}^2+{\sigma^{0,1}}^2\)}\big]}{{1-\E\big[\frac{{\sigma^{0,1}}^2\lambda_{11}}{{\sigma^{1}}^2+{\sigma^{0,1}}^2}\big]}}.$$
	Denote $$A=\frac{\E\big[\frac{{\sigma^{0,1}}^2\lambda_{12}}{{\sigma^{1}}^2+{\sigma^{0,1}}^2}\big]}{1-\E\big[\frac{{\sigma^{0,1}}^2\lambda_{11}}{{\sigma^{1}}^2+{\sigma^{0,1}}^2}\big]},\ B=\frac{\E\big[\frac{\mu^1\sigma^{0,1}}{\delta^1\({\sigma^{1}}^2+{\sigma^{0,1}}^2\)}\big]}{{1-\E\big[\frac{{\sigma^{0,1}}^2\lambda_{11}}{{\sigma^{1}}^2+{\sigma^{0,1}}^2}\big]}}.$$
	Similarly, we find that $\pi^{1,i,*}_{N_1}(u,v)$ is linear with respect to $u$ and independent of $v$.
	Using the first equation of Equation \eqref{eq: equi of 2*N}, we have 
	$x_1=A^Nx_2+B^N$,
	where $$A^N=\frac{\frac{1}{N}\Sigma_{i =1}^{N}\frac{{\sigma^{0,1,i}}^2\lambda^{i}_{1,2}}{(1-\frac{\lambda^{i}_{1,1}}{N}){\sigma^{1,i}}^2+{\sigma^{0,1,i}}^2}}{1-\frac{1}{N}\Sigma_{i =1}^{N}\frac{{\sigma^{0,1,i}}^2\lambda^{i}_{1,1}}{(1-\frac{\lambda^{i}_{1,1}}{N}){\sigma^{1,i}}^2+{\sigma^{0,1,i}}^2}},\ B^N=\frac{\frac{1}{N}\Sigma_{i =1}^{N}\frac{\mu_t^{1,i}\sigma^{0,1,i}}{\delta^{1,i}((1-\frac{\lambda^{i}_{1,1}}{N}){\sigma^{1,i}}^2+{\sigma^{0,1,i}}^2)}}{1-\frac{1}{N}\Sigma_{i =1}^{N}\frac{{\sigma^{0,1,i}}^2\lambda^{i}_{1,1}}{(1-\frac{\lambda^{i}_{1,1}}{N}){\sigma^{1,i}}^2+{\sigma^{0,1,i}}^2}}.$$
	When $\lambda_{12}=0$, we have $A^N=A=0$ and $x^{1,*}_N=B_N,x^{1,*}=B$. The equilibrium strategies in the N-player game and MFG are given by
	$$\begin{aligned}
		\pi^{i,*}_{N,t}
		=&\frac{{\sigma^{0,1,i}}}{(1-\frac{\lambda^{i}_{1,1}}{N_1}){\sigma^{1,i}}^2+{\sigma^{0,1,i}}^2}\lambda^{i}_{1,1}B^N+\frac{\mu^{1,i}}{\delta^{1,i}((1-\frac{\lambda^{i}_{1,1}}{N_1}){\sigma^{1,i}}^2+{\sigma^{0,1,i}}^2)}\\
		=&\frac{{\sigma^{0,1,i}}}{(1-\frac{\lambda^{i}_{1,1}}{N_1}){\sigma^{1,i}}^2+{\sigma^{0,1,i}}^2}\lambda^{i}_{1,1}\frac{\frac{1}{N}\Sigma_{i =1}^{N}\frac{\mu_t^{1,i}\sigma^{0,1,i}}{\delta^{1,i}((1-\frac{\lambda^{i}_{1,1}}{N}){\sigma^{1,i}}^2+{\sigma^{0,1,i}}^2)}}{1-\frac{1}{N}\Sigma_{i =1}^{N}\frac{{\sigma^{0,1,i}}^2\lambda^{i}_{1,1}}{(1-\frac{\lambda^{i}_{1,1}}{N}){\sigma^{1,i}}^2+{\sigma^{0,1,i}}^2}}\\&+\frac{\mu^{1,i}}{\delta^{1,i}((1-\frac{\lambda^{i}_{1,1}}{N_1}){\sigma^{1,i}}^2+{\sigma^{0,1,i}}^2)}
	\end{aligned}
	$$
	and
	\begin{equation}\label{eq: num cir1 MFG eq}
		\begin{aligned}
			\pi^{1,*}_t=&
			\frac{{\sigma^{0,1}}}{{\sigma^1}^2+{\sigma^{0,1}}^2}\lambda_{1,1}B+\frac{\mu^1}{\delta^1({\sigma^1}^2+{\sigma^{0,1}}^2)}\\
			=&\frac{{\sigma^{0,1}}}{{\sigma^1}^2+{\sigma^{0,1}}^2}\lambda_{1,1}\frac{\E[\frac{\mu^1\sigma^{0,1}}{\delta^1({\sigma^{1}}^2+{\sigma^{0,1}}^2)}]}{{1-\E[\frac{{\sigma^{0,1}}^2\lambda_{11}}{{\sigma^{1}}^2+{\sigma^{0,1}}^2}]}}+\frac{\mu^1}{\delta^1({\sigma^1}^2+{\sigma^{0,1}}^2)}.
		\end{aligned}
	\end{equation}
	This result coincides with the main findings in {Lacker and Zariphopoulou \cite{LD2019mean}}, where they used the notations $B=\frac{\varphi}{1-\psi}$ and $B^N=\frac{\varphi_N}{1-\psi_N}$. Within our notation, the meaning of $x^*=B=\frac{\varphi}{1-\psi}$ or $x^*_N=B^N=\frac{\varphi_N}{1-\psi_N}$ is more straightforward: they represent the population mean of the load on Brownian risk. 
	
	Another interesting observation is that the MFE of the population 1 in our multiple population setting is $\pi^{1,*}(u^1_t)=\frac{{\sigma^{0,1}}}{{\sigma^1}^2+{\sigma^{0,1}}^2}u^1_t+\frac{\mu^1}{\delta^1\({\sigma^1}^2+{\sigma^{0,1}}^2\)}$, where $u^1_t=\lambda_{1,1}x_1+\lambda_{1,2}x_2$ representing the weighted mean of risk load on Brownian risk. This implies that when considering the performance of other population with extra Poisson common noise, the population 1 only cares about the Brownian common noise and the idiosyncratic noise within the population-1. The Poisson common noise attached to the population 2 only affects the final value function but not the MFE of the population 1. This is a reasonable result because the population 1 can only invest in the underlying assets drive by the Brownian idiosyncratic and common noise and their own Poisson idiosyncratic noise.

    On the other hand, if we consider $\lambda_{2,1}=0$, which indicates that the population 2 does not concern the performance of first population, but the Poisson common noise is considered. Then the MFE satisfies the simplified characterization. 
	\begin{corollary}
		Let $x^*_2$ and $y^*$ be the solution of  equation system
		\begin{equation}\label{eq: numerical mfg equi p2}
			\left\{\begin{aligned}
				&\E\[\pi^{2,*}(\lambda_{22}x_2,\lambda_{22}y)\sigma^{0,2}\]=x_2,\\
				&\E\[\pi^{2,*}(\lambda_{22}x_2,\lambda_{22}y)\gamma\]=y.\\
			\end{aligned}\right.
		\end{equation}
		Then $\pi^{2,*}(\lambda_{22}x^*_{2},\lambda_{22}y^*)$ is the MFE of population 2.
	\end{corollary}
	Similarly, if we consider $\lambda_{1,2}=0$ that population 1 does not care the performance of the second population, we have the simplified result. 
	\begin{corollary}
		Let $x^*_1$ be the solution of the equation
		\begin{equation}\label{eq: numerical mfg equi p1}
			\E\[\pi^{1,*}(\lambda_{11}x_1,0)\sigma^{0,1}\]=x_1.
		\end{equation}
		Then $\pi^{1,*}(\lambda_{11}x^*_{1},0)$ is the MFE of Population 1.
	\end{corollary}

	\subsection{Numerical examples}
	In this subsection, we present some numerical results to illustrate the different behaviors that arise when jump risk is presented. We assume that all the variables in the type vector follow a truncated Gaussian distribution, truncated in $[mean-3std, mean+3std]$. To make our parameters meaningful, we adopt the estimated parameters  obtained in {Honore \cite{honore1998pitfalls}}. We also assume that two independent Brownian motions have the same volatility, with total volatility equal to that in {Honore \cite{honore1998pitfalls}}. All the parameters are scaled to an annualized approximation. We assume that all the parameters have a standard deviation of 10\% of their means. The following table lists the parameters used.
An interesting fact is that the estimated $\gamma$ is negative.

	\begin{table}[h]
		\centering
		\begin{tabular}{c c c c c c}
			\hline
			$\delta$ & $\gamma$ & $\nu$ & $\sigma$ & $\sigma^0$ & $\mu$ \\
			\hline
			$1.0$ & $-0.04$ & $4.7$ & $0.11$ & $0.11$ & $0.25$ \\
			\hline
		\end{tabular}
		\caption{Mean of Parameters}
		\label{tab:params}
	\end{table}
	  We consider the following scenarios:
	\begin{enumerate}
		\item Population 1 with only Brownian idiosyncratic and common noise. The equilibrium is described by Equation \eqref{eq: num cir1 MFG eq}.
		\item Population 1 with Brownian+Poisson idiosyncratic risk and only Brownian common noise. The equilibrium is described by Equation \eqref{eq: numerical mfg equi p1}.
		\item Population 2  with Brownian idiosyncratic risk as well as Brownian+Poisson common noise. The equilibrium is described by Equation \eqref{eq: numerical mfg equi p2}.
	\end{enumerate}
	We use $\pi^{*,d_1}(u,v),\pi^{*,d_2}(u,v)$ and $\pi^{*,d_3}(u,v)$ to denote the MFE in above three cases.
	
\begin{figure}[h]
\centering
\includegraphics[width=\textwidth]{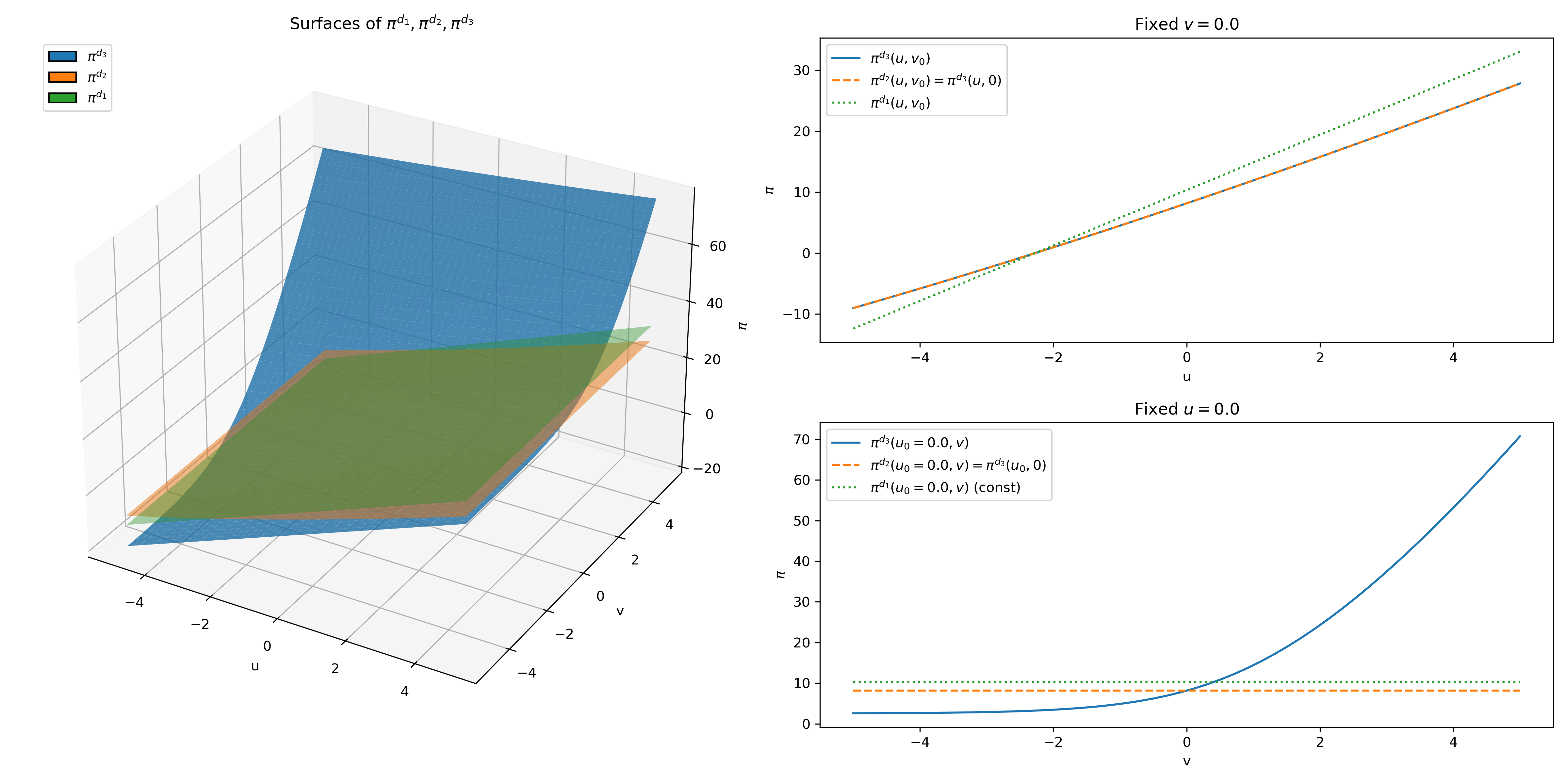}
\caption{Plots of $\pi^*$ with positive $\gamma$.}
\label{fig: pi 2}
\end{figure}

	First, we plot $\pi(u,v)$, with $u$ and $v$ representing the relative mean loadings on the Brownian and Poisson common noises, respectively. As $\pi^*$ depends on several parameters, we only show the plot for type vectors fixed at their mean values. Figure \ref{fig: pi 1} shows the results:
\begin{figure}[h]
\centering
\includegraphics[width=\textwidth]{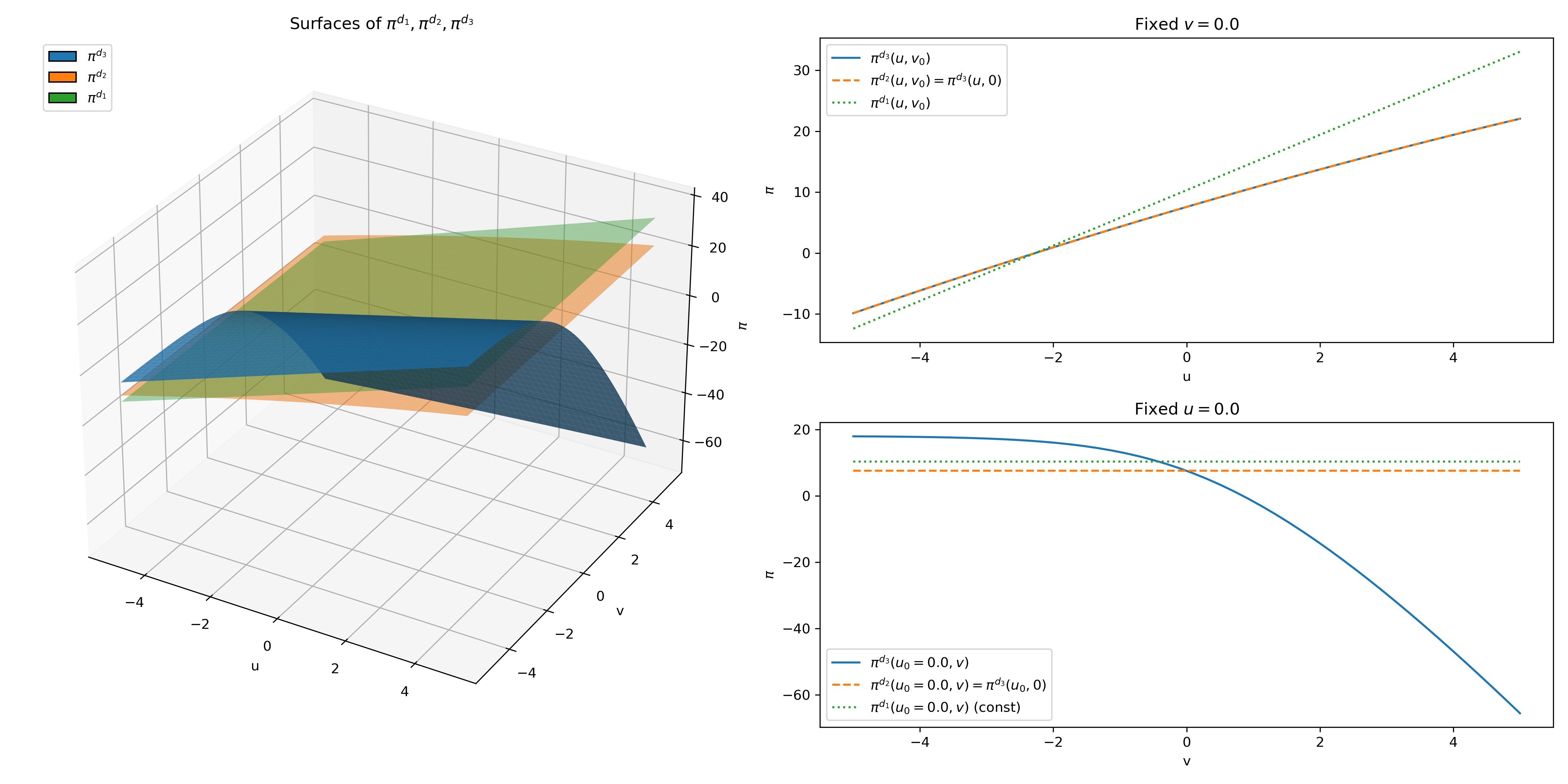}
\caption{Plots of $\pi^*$.}
\label{fig: pi 1}
\end{figure}
The left subplot provides the 3D plot of $\pi^*$ with respect to $(u,v)$, while the right subplots illustrate the influence of $u$ and $v$, with the other parameters fixed. As we can see, both $u$ and $v$ exhibit nonlinear effects on the optimal solution. Specifically, $u$ has a positive influence, while $v$ has a negative influence. Notably, the behavior of the MFE differs substantially when Poisson common noise is considered (the blue surface in the left figure). The influence of $v$ is more nonlinear than that of $u$. The results in the right figures indicate that the greater the Brownian motion exposure in the relative performance process, the higher the Brownian motion risk in the MFE. In contrast, Poisson risk exhibits the opposite effect (recall that $\gamma<0$ here). The green line in the upper-right subplot represents the MFE when $\gamma=0$, i.e., the case with only Brownian motion. In this case, the MFE is linear. We observe that when Poisson risk is introduced, the investor becomes more risk-averse as the relative performance process has a higher loading on Brownian motion. We also show the result with $\gamma=0.04>0$ in Figure \ref{fig: pi 2}.
As we can see, the influence of $v$ becomes positive, while the difference between the linear strategy and the nonlinear strategy with Poisson risk remains similar.
	
\begin{figure}[h]
\centering
\includegraphics[width=\textwidth]{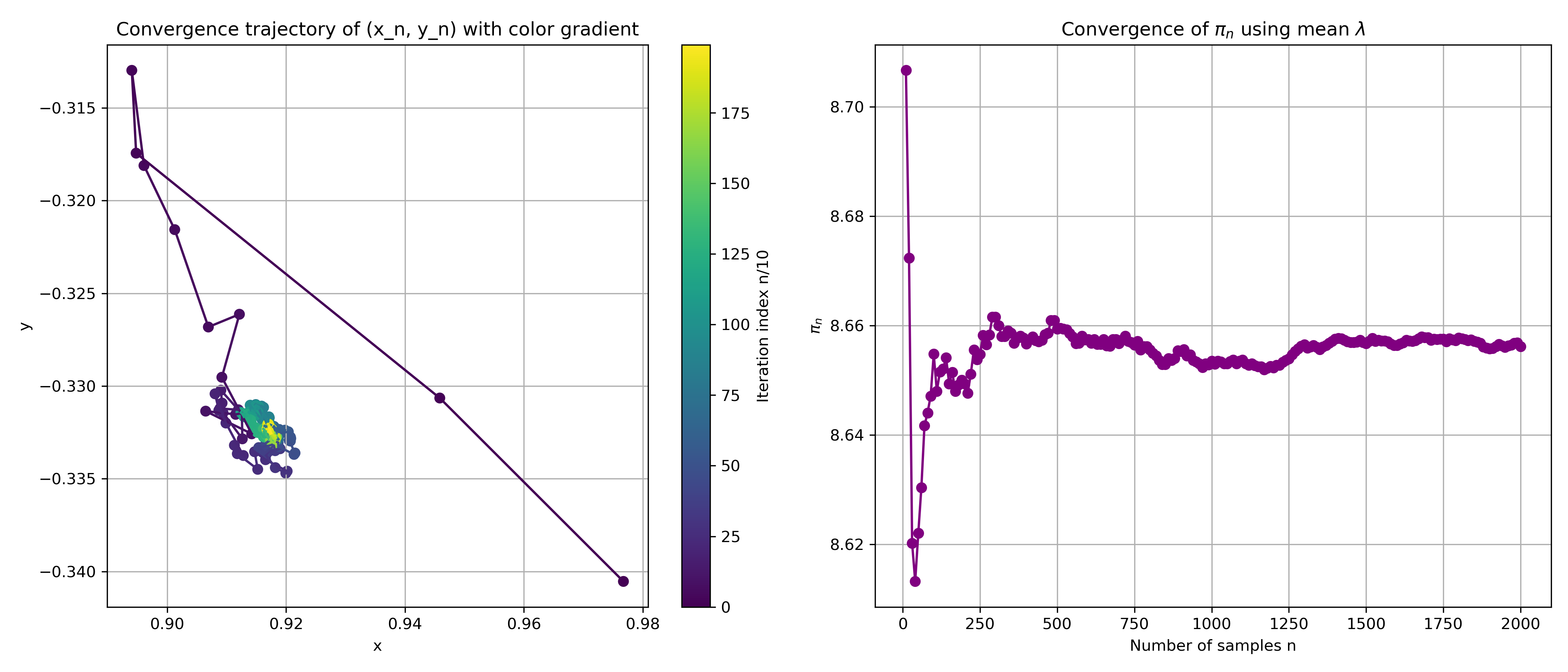}
\caption{Convergence of $(x,y)$ and the equilibrium strategy $\pi^*$.}
\label{fig: converge}
\end{figure}

	Second, we present the equilibrium results of the MFG. Similarly, we assume that $\lambda$ follows a truncated Gaussian distribution with mean $0.2$ and standard deviation $0.02$.
Since it is impossible to obtain the closed-form solution of Equation \eqref{eq: numerical mfg equi p2}, we adopt an approximation approach, which follows exactly the same convergence result proved in Section \ref{sec: convergence}. We randomly sample $N_1$ and $N_2$ different type vectors from the corresponding distributions. By setting $N_1 = N_2 = N = 2000$, we obtain the desired values $x^* = 0.917$ and $y^* = -0.332$ from the approximated equation. The equilibrium strategy $\pi^*$ for an investor with the mean type vector is $\pi^* = 8.656$. We plot the convergence sequences in Figure \ref{fig: converge}.
The left subplot in Figure \ref{fig: converge} illustrates the convergence trajectory of $(x_n, y_n)$, computed using samples $1 \sim n$. The right subplot shows the convergence of $\pi^*$, calculated using the mean of the parameters as the type vector. Meanwhile, we find that the MFE $\pi^*$ without relative performance ($\gamma = 0$) is $7.58$, which is approximately 12\% lower than the previous value of $\pi^* = 8.656$. This result indicates that competition within the population significantly increases the amount of wealth invested in risky assets.

\begin{figure}[h]
\centering
\includegraphics[width=0.6\textwidth]{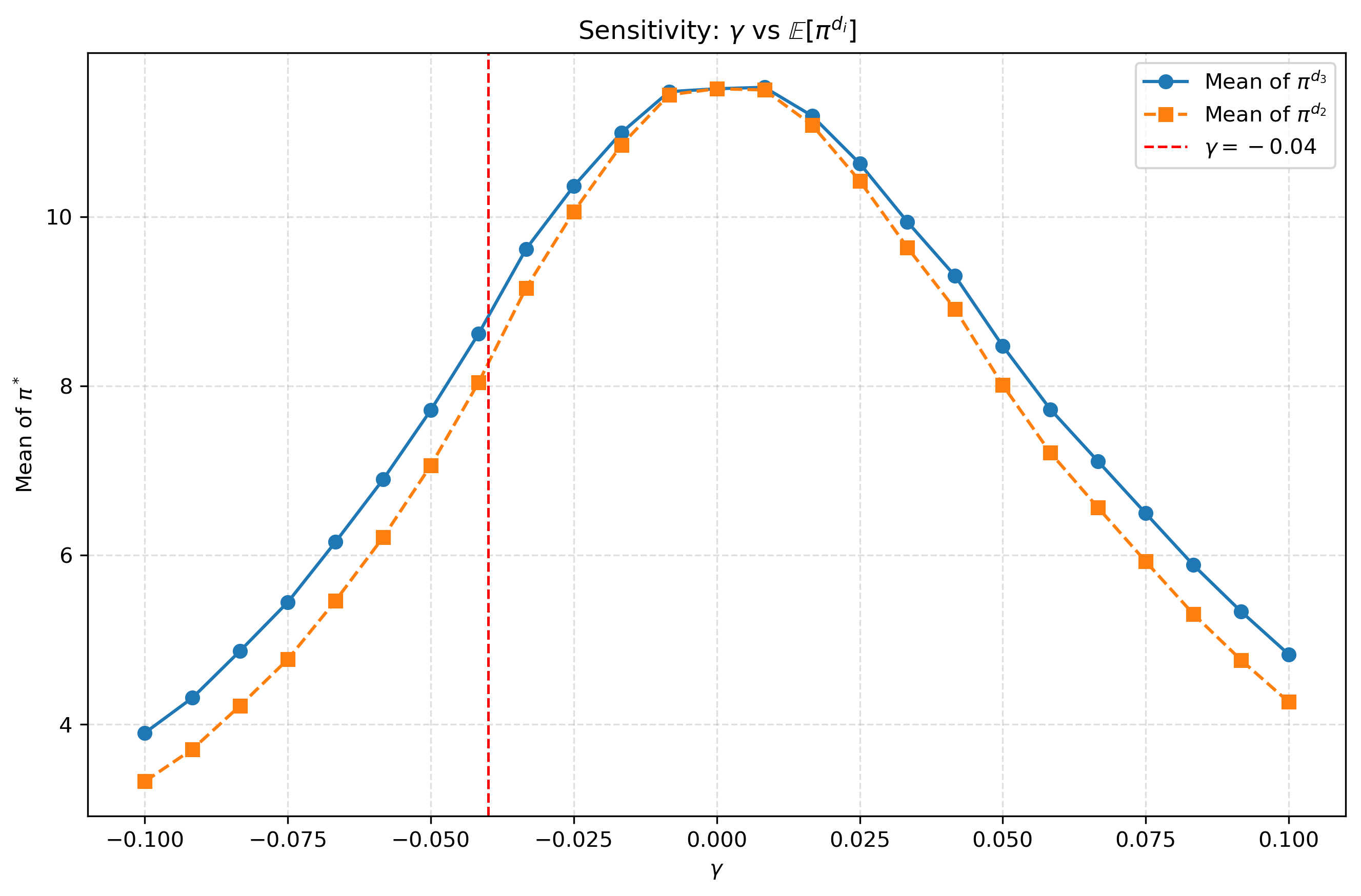}
\caption{Sensitivity with respect to $\gamma$}
\label{fig: gamma vs pi}
\end{figure}

    Next, we present the sensitivity analysis of the MFE with respect to different model parameters. Because the equilibrium strategy is a random variable which is hard to compare, we compare parameters with the expectation of the equilibrium strategy. In Figure \ref{fig: gamma vs pi}, we vary the mean of $\gamma$ (std remains 10\% of the mean), and conduct experiments to show the change in the expectation of the MFE $\pi^{d_2}$ and $\pi^{d_3}$. The results show that when $|\gamma|$ is larger, the investors intend to invest less in the risky asset (due to the increased risk). The curve of $\pi^{d_3}$ is always above $\pi^{d_2}$, which indicates that the Poisson common noise has less adverse effect than the Poisson idiosyncratic jump risk in the MFE allocation in the risky asset. In the presence of idiosyncratic jump risk, the investor will be more conservative in the portfolio decision at the equilibrium state.

\begin{figure}[h]
\centering
\includegraphics[width=0.6\textwidth]{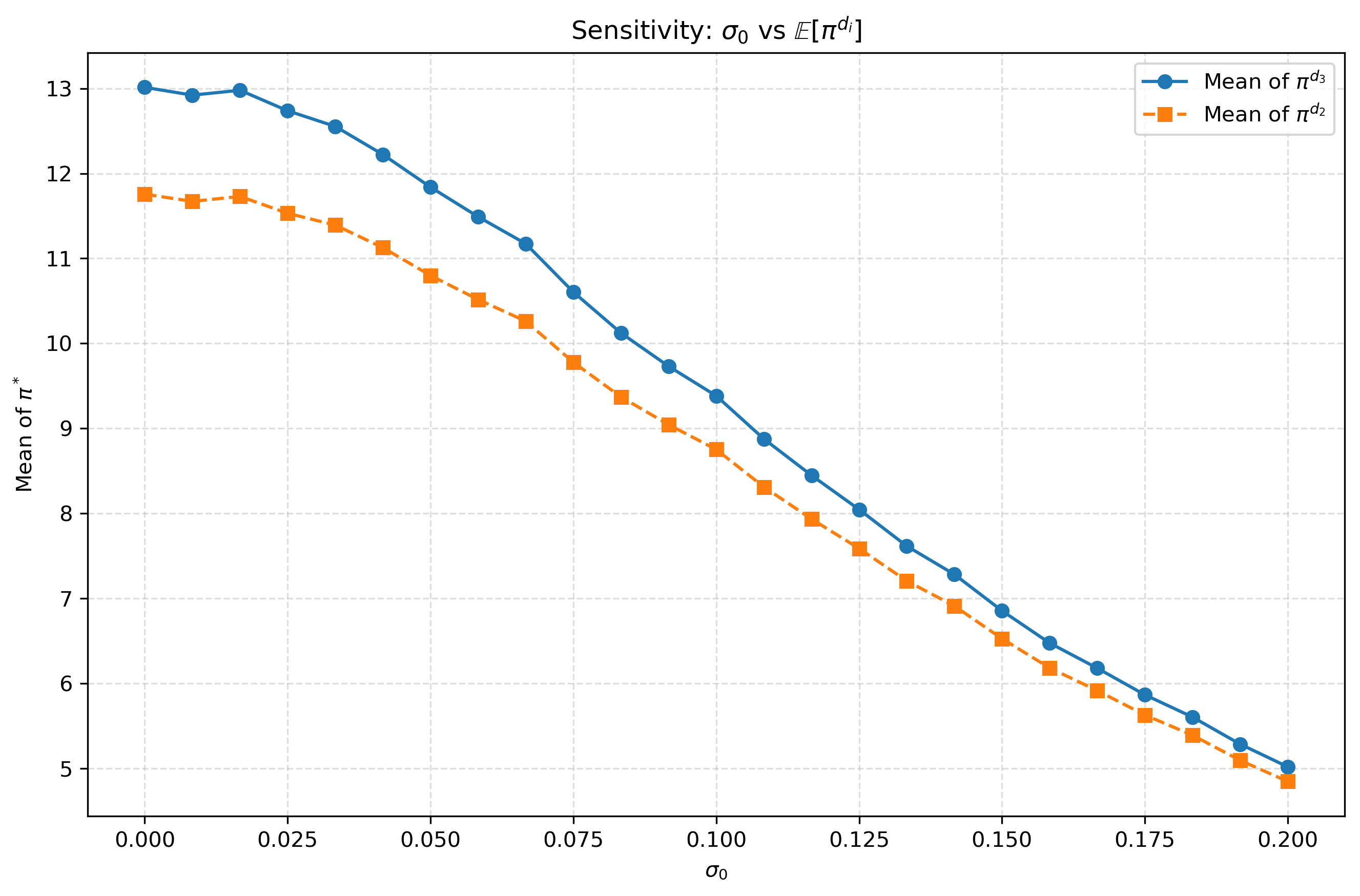}
\caption{Sensitivity with respect to $\sigma^0$}
\label{fig: sigma0 vs pi}
\end{figure}

\vskip 1cm
In Figure \ref{fig: sigma0 vs pi}, the influence of $\sigma^0$ is considered. As we can see, when $\sigma^0$ increases, the MFE strategy decreases, which is similar to the influence of $\gamma$. Also, we can see that the competition encourages investors to allocate more in the risky asset.

\begin{figure}[h]
\centering
\includegraphics[width=0.6\textwidth]{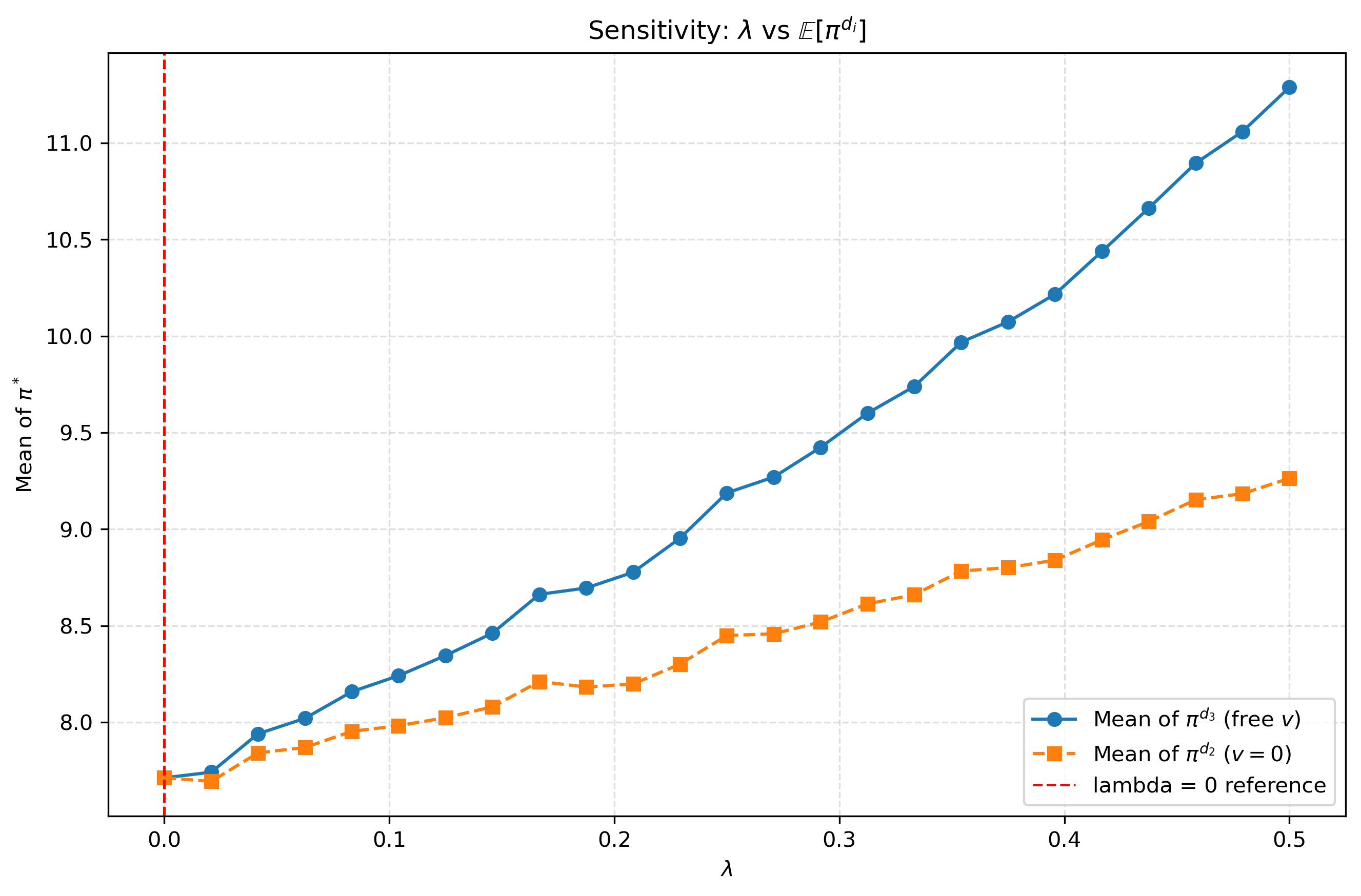}
\caption{Sensitivity with respect to $\lambda$}
\label{fig: lambda vs pi}
\end{figure}

Finally, we show the influence of $\lambda$ in Figure \ref{fig: lambda vs pi}. Clearly, investors intend to allocate more in the risky asset when $\lambda$ is high, which implies that the level of relative performance has a positive influence on the MFE.
	\section*{Acknowledgments.}
	Zongxia Liang is supported by the National Natural Science Foundation of China under grant no.\ 12271290.
	Xiang Yu is supported by the Hong Kong RGC General Research Fund (GRF) under grant no.\ 15211524, 
	the Hong Kong Polytechnic University research grant under no.\ P0045654, 
	and the Research Centre for Quantitative Finance at the Hong Kong Polytechnic University under grant no.\ P0042708.
\bibliographystyle{plain}
\bibliography{MFGwithjump.bib}

\end{document}